\documentstyle[12pt,amssymb]{article}

\newlength{\minitwocolumn}
\newcommand{\Z}{{\Bbb Z}} 
\newcommand{\C}{{\Bbb C}} 

\newcommand{\F}{{\cal F}}

\newcommand{\cR}{{\cal R}}

\newcommand{\bR}{{\bar{R}}}
\newcommand{\hL}{\widehat{L}}

\renewcommand{\H}{{\cal H}}
\newcommand{\la}{\lambda}

\newcommand{\nn}{{\nonumber}}
\newcommand{\eqref}[1]{(\ref{#1})}
\newcommand{\bea}{\begin{eqnarray}}
\newcommand{\ena}{\end{eqnarray}}
\newcommand{\beit}{\begin{itemize}}
\newcommand{\enit}{\end{itemize}}

\newcommand{\be}{\begin{eqnarray*}}
\newcommand{\en}{\end{eqnarray*}}
\newcommand{\lb}[1]{\label{#1}}

\newcommand{\id}{\hbox{id}}
\newcommand{\Ad}{\mbox{Ad}}
\newcommand{\Diag}{\mbox{Diag}}

\def\infq4p#1{{(#1;q^4,p)_\infty}}

\newcommand{\al}{\alpha}

\newcommand{\vep}{\varepsilon}

\newcommand{\ba}{\bar{\alpha}}

\newcommand{\rv}{{\vee}}
\newcommand{\alv}{\al^{\vee}}

\font\teneufm=eufm10
\font\seveneufm=eufm7
\font\fiveeufm=eufm5
\newfam\eufmfam
\textfont\eufmfam=\teneufm
\scriptfont\eufmfam=\seveneufm
\scriptscriptfont\eufmfam=\fiveeufm

\let\goth\frak
\newcommand{\slth}{\widehat{\goth{sl}}_2}
\newcommand{\slt}{\goth{sl}_2}
\newcommand{\slnh}{\widehat{\goth{sl}}_N}

\newcommand{\g}{\goth{g}}

\newcommand{\Aqp}{{\cal A}_{q,p}}

\newcommand{\Bqla}{{{\cal B}_{q,\lambda}}}

\newcommand{\att}{A^{(2)}_2}

\newcommand{\h}{\goth{h}}

\makeatletter
\@addtoreset{equation}{section}
\makeatother

\newtheorem{thm}{Theorem}[section]
\newtheorem{prop}[thm]{Proposition}

\newtheorem{cor}[thm]{Corollary}

\newtheorem{df}{Definition}[section]
\newtheorem{dfn}[thm]{Definition}

\begin{document}

\

\vspace{0.7cm}
\begin{center}
{\Large{\bf
The Drinfeld Realization of
\\
the Elliptic Quantum Group ${\cal B}_{q,\lambda}
(A_2^{(2)})$
}}

\vspace{1.2cm}
{ Takeo KOJIMA$~^{*}$ and~ Hitoshi KONNO$~^{**}$}

\vspace{0.7cm}
{\it
~*
Department of Mathematics,
College of Science and Technology,\\
Nihon University, Chiyoda-ku, Tokyo
101-0062, Japan.\\
E-mail:kojima@math.cst.nihon-u.ac.jp\\~\\
~**
Department of Mathematics,
Faculty of Integrated Arts and Sciences,\\
Hiroshima University,
Higashi-Hiroshima 739-8521, Japan.
\\
E-mail:konno@mis.hiroshima-u.ac.jp
}

 
\end{center}

\vspace{0.8cm}
\begin{abstract}
We construct a realization of the $L$-operator satisfying the 
$RLL$-relation of the face type elliptic quantum group 
$\Bqla(A_2^{(2)})$. The construction is based on the elliptic analogue 
of the Drinfeld currents of $U_q(\att)$, which forms the elliptic algebra 
$U_{q,p}(A_2^{(2)})$. We give a realization of the elliptic currents
$E(z), F(z)$ and $K(z)$ as a tensor product of the Drinfeld 
currents of $U_q(A_2^{(2)})$ and a Heisenberg algebra. In the level-one 
representation, we also give a free field realization of the elliptic currents. Applying these results, we derive a free field realization of the 
$U_{q,p}(A_2^{(2)})$-analogue of the $\Bqla(A_2^{(2)})$-intertwining 
operators. The resultant operators coincide with those of the vertex 
operators in the dilute $A_L$ model, which is known to be 
a RSOS restriction of the $\att$ face model.  
\end{abstract}

\newpage

\section{Introduction}

An elliptic quantum groups is a quasi-triangular quasi-Hopf 
algebra obtained as a quasi-Hopf deformation of the affine quantum 
group $U_q(\g)$ by the twistor satisfying the shifted cocycle condition
\cite{EF,Fronsdal,JKOS1}. It is conjectured in \cite{FIJKMY,JKOS1} that the 
representation theory of the elliptic quantum groups of both the vertex type 
$\Aqp(\slnh)$ and the face type $\Bqla(\g)$, $\g$ being an affine Lie algebra, 
enables us to perform an algebraic analysis of the corresponding two 
dimensional solvable lattice models in the sense of Jimbo and Miwa \cite{JM}. 
In order to perform 
the analysis, we need to construct explicit representations of both finite 
and infinite dimensional. For this purpose, the Drinfeld realization of the 
quantum groups is known to provide a relevant framework. In the previous 
papers\cite{Konno,JKOS2,KoKo}, we constructed the Drinfeld realization of the 
face type elliptic quantum group $\Bqla(\slnh)$ based on the elliptic algebra
$U_{q,p}(\slnh)$. The Drinfeld generators have both finite and infinite 
dimensional representations suitable for the calculation of the correlation 
functions. 

In this paper, we investigate the same problem for $\Bqla(A^{(2)}_2)$, 
the face type elliptic 
quantum group associated with the twisted affine 
Lie algebra $A^{(2)}_2$. 
We first construct the elliptic algebra 
$U_{q,p}(A^{(2)}_2)$ as the algebra of the elliptic analogue of the Drinfeld
currents of $U_q(A^{(2)}_2)$. 
Basically, the idea given in Appendix A of \cite{JKOS2} can be 
applied to the twisted case.
Namely, dressing the Drinfeld currents of $U_q(A^{(2)}_2)$ by the bosons 
$a_m\ (m\in\Z_{\not=0})$ in $U_q(A^{(2)}_2)$ 
and taking a tensor product with a certain 
Heisenberg algebra $\C\{\H\}$ generated by $P, Q$, 
which commutes with $U_q(\att)$, we obtain
the elliptic Drinfeld currents. However, we formulate the elliptic 
algebra $U_{q,p}(A^{(2)}_2)$ in an extended form by introducing the new 
currents $K(u)$, which enables the $RLL$-formulation of $U_{q,p}(\att)$.
Then we discuss a connection  
between $U_{q,p}(A^{(2)}_2)$ and $\Bqla(\att)$. 
We derive the dynamical $RLL$-relation of $\Bqla(A^{(2)}_2)$ from 
the $RLL$-relation of $U_{q,p}(\att)$ 
by removing a half of the generator $Q$ of the Heisenberg algebra
 and identifying $P$ with the dynamical parameter 
 in $\Bqla(\att)$. We hence find a structure
of $U_{q,p}(A^{(2)}_2)$ roughly given by 
$``\Bqla(A^{(2)}_2)\otimes \C\{\H\}"$, 
and in this sense, we regard $U_{q,p}(A^{(2)}_2)$ as the Drinfeld realization 
of $\Bqla(\att)$.

Although the above tensor structure does not preserve the 
coalgebra structure of $\Bqla(A^{(2)}_2)$, the same 
tensor structure enables us to convert the algebraic objects 
of $\Bqla(A^{(2)}_2)$, such as the intertwining operators, to the 
$U_{q,p}(A^{(2)}_2)$ counterparts. 
In the known cases, it is true that the $U_{q,p}(\g)$ 
counterparts of the $\Bqla(\g)$ intertwining 
operators play the role of vertex operators in the restricted solid 
on solid (RSOS) model associated with $\g$. We call such ``intertwining" 
operator of $U_{q,p}(\g)$ the vertex operator of $U_{q,p}(\g)$. 
Moreover, the elliptic Drinfeld currents in $U_{q,p}(A^{(2)}_2)$ 
admits a free field realization, which is an elliptic extension of those of 
$U_q(\att)$ obtained in \cite{Matsuno,JingMisra,Hou}. By 
using such realization 
and applying the tensoring procedure, we derive a free field realization of 
the vertex operators of $U_{q,p}(\att)$.

The face model associated with the twisted affine Lie algebra $A^{(2)}_2$ 
was formulated in \cite{Kuniba}. Its RSOS restriction 
is known to be the dilute $A_L$ model\cite{WNS92,WPSN94}. 
The free filed formulation of the dilute $A_L$ model was carried out 
in \cite{HJKOS}. There, however, the construction of the vertex operators 
was done by brute force based on the commutation relations among the vertex 
operators and on a partial result on the elliptic Drinfeld currents.
We here derive the same vertex operators by using the 
the representation theory of $\Bqla(\att)$ and the Drinfeld realization 
given by $U_{q,p}(\att)$.

This paper is organized as follows. 
In the next section, we give a summary of the basic facts on 
the face type elliptic quantum group $\Bqla(A^{(2)}_2)$. 
In Section 3, we present a definition and a realization of the 
elliptic algebra $U_{q,p}(A^{(2)}_2)$. New currents $K(u)$ 
are introduced there. 
In Section 4, we introduce a set of half currents defined 
from the elliptic currents in $U_{q,p}(A^{(2)}_2)$ 
and derive their commutation relations. 
Section 5 is devoted to a construction of a $L$-operator 
 and the $RLL$-formulation of $U_{q,p}(\att)$. 
We then derive the dynamical $RLL$-relation of $\Bqla(A^{(2)}_2)$ 
from $U_{q,p}(\att)$. 
According to this result, in Section 6, 
we discuss a free field realization of the two types of 
vertex operators of the level one $U_{q,p}(A^{(2)}_2)$-modules. 
The final section is devoted to a discussions on some remaining 
problems.
In addition, we have three appendices. 
In Appendix A, we give a summary of the $3$ dimensional   
evaluation representation of $U_{q,p}(A^{(2)}_2)$. 
In Appendix B, we discuss 
the difference equation for the twistor and give a partial results on the 
solutions. Finally, in Appendix C, we give a 
proof of some formulae of commutation relations of the half currents.
 
\newpage
\section{The Elliptic Quantum Group 
${\cal B}_{q,\lambda}(A_2^{(2)})$}

In this section, we summarize some basic facts on the 
face type elliptic quantum group $\Bqla(A_2^{(2)})$ based 
on the results in \cite{JKOS1}. 

\subsection{Notations}
Through this article, we fix a complex number $q\neq0, 0<q<1$ and $p$ given by
\begin{eqnarray*}
p=q^{2r},~~p^*=pq^{-2c}=q^{2r^*}~~(r^*=r-c;~ r,r^* \in 
{\mathbb{R}}_{>0}).
\end{eqnarray*}
We parametrize $p$ as follows.
\begin{eqnarray*}
&&p=e^{-2\pi i/\tau},~p^*=e^{-2\pi i/{\tau^*}}~~(r\tau=r^*\tau^*),\\
&&z=q^{2u}=e^{-2\pi i u/r \tau}
.
\end{eqnarray*}
We often use the following Jacobi theta functions.
\begin{eqnarray*}
&&~[u]=q^{\frac{u^2}{r}-u}{\Theta_p(q^{2u})}=e^{-\frac{\pi i}{4}}\tau^{\frac{1}{2}}q^{-\frac{r}{4}}\ \vartheta_1\left(\frac{u}{r}\biggl|\tau\right),\\
&&~[u]_+=q^{\frac{u^2}{r}-u}
{\Theta_p(-q^{2u})}=e^{-\frac{\pi i}{4}}\tau^{\frac{1}{2}}q^{-\frac{r}{4}}\ 
\vartheta_0\left(\frac{u}{r}\biggl|\tau\right),
\end{eqnarray*}
$[u]^*=[u]|_{r\to r^*, \tau\to\tau^*}$ and 
$[u]_+^*=[u]_+|_{r\to r^*, \tau\to\tau^*}$.
Here 
\begin{eqnarray*}
&&\Theta_p(z)=(z,p)_\infty (pz^{-1};p)_\infty
(p;p)_\infty,\\
&&(z;t_1,\cdots,t_k)_\infty=
\prod_{n_1,\cdots,n_k \geq 0}(1-zt_1^{n_1}\cdots t_k^{n_k}).
\end{eqnarray*}
The theta functions satisfy $[-u]=-[u], [-u]_+=[u]_+$ 
and the quasi-periodicity property
\begin{eqnarray}
&&~[u+r]=-[u],~~[u+r\tau]=-e^{-\pi i \tau-\frac{2\pi i u}{r}}[u],\\
&&
~[u+r]_+=[u]_+,~~[u+r\tau]_+=e^{-\pi i \tau-\frac{2\pi i u}{r}}
[u]_+,\\
&&[u+\frac{r\tau}{2}]=ie^{-\pi i(u/r+\tau/4)}[u]_+.
\end{eqnarray}
We use the following normalization for the contour integration.  
\begin{eqnarray}
\oint_{C_0}\frac{dz}{2\pi i z}\frac{1}{[-u]}=1,\qquad 
\oint_{C_0}\frac{dz}{2\pi i z}\frac{1}{[-u]^*}=\frac{[u]}{[u]^*}
\Bigl\vert_{u\to 0}
\end{eqnarray}
where $C_0$ is a simple closed curve 
in the $u$-plane encircling
$u=0$ anticlockwise.

\subsection
{Definition of the elliptic quantum group $\Bqla(A_2^{(2)})$}

Let $U_{q}(A_2^{(2)})$ 
be the standard affine quantum group,
associated with the Cartan matrix 
\begin{eqnarray}
A=\left(\begin{array}{cc}
2&-1\\
-4&2
\end{array}\right).
\end{eqnarray}
The label of $A$ is $a_0=1, a_1=2$ and colabel is 
$a_0^\rv=2, a_1^\rv=1$. 
Let $B=(b_{ij})$ be the symmetrized Cartan matrix 
$b_{ij}=\frac{a_i^\rv}{a_i}a_{ij}$.  
We identify $\h=\C \al^{\vee}_0\oplus \C \al^\rv_1\oplus \C d$ 
and $\h^*=\C \al_0\oplus \C \al_1\oplus \C \Lambda_0$ via 
the standard invariant bilinear form $(\ ,\ )$ 
given  on $\h$ and $\h^*$ as follows.
\be
&&(\al^\rv_i ,\al^\rv_j )=a_{ij}\frac{a_j}{a_j^\rv} \quad (0\leq i,j\leq 1),\\
&&(\al^\rv_i,d)=\delta_{i,0},\qquad (d,d)=0,\\
&&(\al_i ,\al_j )=\frac{a_i^\rv}{a_i}a_{ij} \quad (0\leq i,j\leq 1),\\
&&(\al_i,\Lambda_0)=\delta_{i,0},\qquad (\Lambda_0,\Lambda_0)=0,\\
\en
The central element is given by $c=2\al_0^\rv+\al^\rv_1$.
Let us set $\delta=\al_0+2\al_1$. Then the following  
relations hold.
\be
&&(\delta,d)=1,\qquad (\delta,\delta)=0,\qquad (c,d)=2,\qquad (c,c)=0.
\en
The identification between $\h$ and $\h^*$ is given explicitly by  
$\al^\rv_i=\frac{2\al_i}{(\al_i,\al_i)}$, $c=\delta$ and $d=2\Lambda_0$. 
Under this, we use $\{\hat{h}_l\}_{l=1,2,3}=\{d, c, \al_1^\rv \}$ 
as a basis of $\h$ and 
its dual basis $\{\hat{h}^l\}_{l=1,2,3}=\{c/2, d/2, \al_1^\rv/2 \}$.  
Our conventions of coalgebra structure of 
$U_q(A^{(2)}_2)$ 
follows \cite{JKOS1}.  The coproduct, counit, antipode are denoted by 
$\Delta$,   $\vep$ and $S$, respectively.

The face type elliptic quantum group 
$\Bqla(A_2^{(2)})$ 
is a quasi-triangular quasi-Hopf algebra obtained from  
$U_{q}(A_2^{(2)})$ 
by the deformation via the face type twistor $F(\la)\ (\la \in \h)$.
The twistor $F(\la)$ is an invertible element in $U_q(\att)\otimes U_q(\att)$ 
satisfying 
\bea
&&({\rm id}\otimes \vep )F(\la)=1=F(\la)(\vep \otimes {\rm id}),\\
&&F^{(12)}(\la)(\Delta\otimes {\rm id})F(\la)
=F^{(23)}(\la+h^{(1)})({\rm id} \otimes \Delta)F(\la).
\lb{facecocy}
\lb{epF}
\ena
where $\la=\sum_l\la_l\hat{h}^l\  (\la_l\in \C)$, 
$\la+h^{(1)}=\sum_l(\la_l+\hat{h}_l^{(1)})\hat{h}^l$ and 
$\hat{h}_l^{(1)}=\hat{h}_l\otimes 1\otimes 1$.  
An explicit construction of the twistor $F(\la)$ is given in \cite{JKOS1}.
As an associative algebra,  
$\Bqla(A_2^{(2)})$ is isomorphic to $U_q(A_2^{(2)})$, 
but the coalgebra structure is deformed in the following way.
\bea
&&\Delta_\la(x)=
F(\la)\Delta(x)F(\la)^{-1}\qquad \forall x \in U_q(A_2^{(2)}).
\ena
$\Delta_\la$ satisfies a weaker coassociativity
\bea
&&({\rm id}\otimes \Delta_\la )\Delta_\la(x)=\Phi(\la)
(\Delta_\la \otimes {\rm id})\Delta_\la(x)\Phi(\la)^{-1}
\qquad \forall\ x \in U_q(A_2^{(2)}),\\
&&\Phi(\la)=F^{(23)}(\la)F^{(23)}(\la+h^{(1)})^{-1}.
\ena
Let $\cR$ be the universal $R$ matrix of $U_q(A^{(2)}_2)$. 
The universal $R$ matrix of $\Bqla(A^{(2)}_2)$ is given by 
\bea
&&\cR(\la)=F^{(21)}(\la)\cR F^{(12)}(\la)^{-1}.\lb{frf}
\ena

\begin{dfn}{\bf (Elliptic quantum group $\Bqla(A_2^{(2)})$)}~~
The face type elliptic quantum group $\Bqla(A_2^{(2)})$ 
is a quasi-triangular
quasi-Hopf algebra 
$(\Bqla(A_2^{(2)}),$
$~\Delta_\la, \vep,~S,~\Phi(\la),~\alpha,~\beta,~\cR(\la))$,
where $\alpha,\ \beta$ are defined by 
\bea
&&\alpha=\sum_iS(k_i)l_i,\quad \beta=\sum_i m_iS(n_i).
\ena
Here we set $\sum_ik_i\otimes l_i=F(\la)^{-1},\ 
\sum_im_i\otimes n_i=F(\la)$.
\end{dfn}

The universal $R$ 
matrix $\cR(\la)$ satisfies the dynamical Yang-Baxter equation.
\begin{equation}
\cR^{(12)}(\la+ h^{(3)} )\cR^{(13)}(\la)\cR^{(23)}(\la+h^{(1)})
=\cR^{(23)}(\la)\cR^{(13)}(\la+h^{(2)})\cR^{(12)}(\la).
\label{DYBE}
\end{equation}
Let $(\pi_{V,z},V_z),\ V_z=V\otimes \C[z,z^{-1}]$ 
be a (finite dimensional) 
evaluation representation 
of $U_q$. Taking images of $\cR$, we define a $R$-matrix 
$R^+_{VW}(z,\la)$ and a $L$-operator $L^+_V(z,\la)$ as follows.
\bea
&&R^{+}_{VW}(z_1/z_2,\la)=
\left(\pi_{V,z_1}\otimes\pi_{W,z_2}\right)
q^{c\otimes d+d\otimes c}
\cR(\la),\\
&&L_V^+(z,\la)=
\left(\pi_{V,z}\otimes {\rm id} \right)
q^{c\otimes d+d\otimes c}\cR(\la).
\ena
Then from \eqref{DYBE}, we have
 the following dynamical $RLL$-relation. 
\bea
&&R^{+}_{VW}(z_1/z_2,\la+ h)
L_V^{+}(z_1,\la)L_W^{+}(z_2,\la+ h^{(1)})
=
L_W^{+}(z_2,\la)L_V^{+}(z_1,\la+ h^{(2)})
R^{+}_{VW}(z_1/z_2,\la).\nn\\
\lb{DRLL}
\ena
Note that in $\Bqla(A_2^{(2)})$, 
$L^+_V(z,\la)$ and $L^-_V(z,\la)=
\left(\pi_{V,z}\otimes {\rm id} \right)
\cR^{(21)}(\la)^{-1}q^{-c\otimes d-d\otimes c}$ 
are not independent 
operators (Proposition 4.3 in \cite{JKOS1}). 
Hence one 
dynamical $RLL$-relation \eqref{DRLL} characterizes the algebra 
$\Bqla(A_2^{(2)})$ 
completely in the sense of Reshetikhin and Semenov-Tian-Shansky
\cite{RS}.

Through this paper, we parametrize the dynamical variable $\lambda$ as
\bea 
&&\lambda=(r^*+3)d+s'c+\frac{1}{2}\left(s+\frac{r\tau}{2}\right)\al_1^\rv
\qquad 
(r^*\equiv r-c).\lb{sj}
\ena 
Under this, we set  $F(r^*,s)\equiv F(\la)$ and 
$\cR(r^*,s) \equiv \cR(\la)$.
Since $c$ is central, no $s'$ dependence should appear. 
The dynamical shift $\lambda \to \lambda+h$  
with $h= cd+(\al_1^{\rv})^2/2$, 
changes  the universal $R$-matrix $\cR(r^*,s)$ 
to 
$\cR(r,s+\al_1^{\rv})$.
Let us take $(\pi_{V,z}, V_z)$ to 
be the  evaluation representation associated with the 
vector representation $V\cong \C^3$ of $U_q(A_2^{(2)})$ 
(see Appendix \ref{Evaluation}).
We set 
\be
&&R^+(u,s+\alv_1)=(\pi_{V,z_1}\otimes \pi_{V,z_2})
q^{c\otimes d+d\otimes c}
\cR(r,s+\alv_1),\\
&&L^+(u,s)=
(\pi_{V,z}\otimes \id)
q^{c\otimes d+d\otimes c}\cR(r^*,s),
\en
where $z_1/z_2=q^{2u}$, $u=u_1-u_2$. 
From \eqref{frf}, we can obtain an explicit expression of $R^+(u,s)$, 
if we know the finite dimensional representation of the twistor 
$(\pi_{V,z_1}\otimes \pi_{V,z_2})
F(r,s)$.  In principle, one can obtain such representation by solving the  
$q$-difference equation for the twistor \cite{JKOS1}, which is   
similar to the $q$-KZ equation for corresponding $U_q(\g)$.
In the present case, the $q$-difference equation splits into the three parts; 
two $2\times 2$ matrix parts  
and one $3\times 3$ matrix part ( see 
Appendix \ref{Twistor}).
Each $2\times 2$ matrix parts turns out to be 
the same as the one of the twistor for  $\Bqla(A_1^{(1)})$ in the vector 
representation after adjusting some $q$-shift and sign factor, 
whereas we have no known 
solutions for the $3\times 3$ matrix part. Writing down the solutions of 
the $2\times 2$ matrix parts under the parametrization of $\lambda$ (\ref{sj}),  we obtain from \eqref{frf} the corresponding matrix elements of $R^+(u,s)$ 
which  coincide with the corresponding matrix elements of the Boltzmann weight 
 for the $A_2^{(2)}$ face model\cite{Kuniba}.  For 
 the remaining $3\times 3$ matrix part, we  
 conjecture that the same coincidence should occur.   
We hence assume that
the $R$-matrix $R^{+}(u,s)$
is given by the following formula.
\begin{eqnarray}
R^+(u,s)=\rho^+(u)\bar{R}(u,s),
\label{def:R}
\end{eqnarray}
where
\begin{eqnarray}
&&\bar{R}(u,s)=\left(\begin{array}{ccccccccc}
1&0&0&0&0&0&0&0&0\\
0&R_{+0}^{+0}&0&R_{+0}^{0+}&0&0&0&0&0\\
0&0&R_{+-}^{+-}&0&R_{+-}^{00}&0&R_{+-}^{-+}&0&0\\
0&R_{0+}^{+0}&0&R_{0+}^{0+}&0&0&0&0&0\\
0&0&R_{00}^{+-}&0&R_{00}^{00}&0&R_{00}^{-+}&0&0\\
0&0&0&0&0&R_{0-}^{0-}&0&R_{0-}^{-0}&0\\
0&0&R_{-+}^{+-}&0&R_{-+}^{00}&0&R_{-+}^{-+}&0&0\\
0&0&0&0&0&R_{-0}^{0-}&0&R_{-0}^{-0}&0\\
0&0&0&0&0&0&0&0&1
\end{array}\right),\lb{rmat}
\end{eqnarray}
\begin{eqnarray*}
&&R_{+0}^{+0}(u,s)=
-\frac{[s+3/2]_+[s-1/2]_+}{[s+1/2]_+^2}
\frac{[u]}{[u+1]},\\
&&R_{+0}^{0+}(u,s)=
\frac{[s+1/2+u]_+[1]}{[s+1/2]_+[1+u]},\\
&&R_{0+}^{+0}(u,s)=
\frac{[-s-1/2+u]_+[1]}{[-s-1/2]_+[1+u]},\\
&&R^{0+}_{0+}(u,s)=R^{-0}_{-0}(u,s)=
-\frac{[u]}{[1+u]},\\
&&R_{0-}^{0-}(u,s)=
-\frac{[s+1/2]_+[s-3/2]_+}{[s-1/2]_+^2}
\frac{[u]}{[u+1]},\\
&&R_{0-}^{-0}(u,s)=
\frac{[s-1/2+u]_+[1]}{[s-1/2]_+[1+u]},\\
&&R_{-0}^{0-}(u,s)=
\frac{[-s+1/2+u]_+[1]}{[-s+1/2]_+[1+u]},\\
&&R_{+-}^{+-}(u,s)=G_s^+G_s^-
\frac{[1/2+u][u]}{[3/2+u][1+u]},\\
&&R_{+-}^{00}(u,s)=
-G_s^-\frac{[s+1/2]_+[-s-1-u]_+[1][u]}
{[-s+1/2]_+^2[1+u][u+3/2]},\\
&&R_{+-}^{-+}(u,s)=
\frac{[-2s+1-u][1]}{[-2s+1][1+u]}-G_s^-
\frac{[-2s-1/2-u][u][1]}{[-2s+1][3/2+u][1+u]},\\
&&R_{00}^{-+}(u,s)=
-\frac{[-s-1-u]_+[1][u]}
{[s+1/2]_+[1+u][u+3/2]},\\
&&R_{-+}^{-+}(u,s)=
\frac{[1/2+u][u]}
{[3/2+u][1+u]},\\
&&R_{-+}^{00}(u,s)=
-\frac{[s-1-u]_+[1][u]}
{[-s+1/2]_+[1+u][u+3/2]},\\
&&R_{-+}^{+-}(u,s)=
\frac{[2s+1-u][1]}{[2s+1][1+u]}-G_s^+
\frac{[2s-1/2-u][u][1]}{[2s+1][3/2+u][1+u]}
,\\
&&R_{00}^{+-}(u,s)=
-G_s^+ \frac{[-s+1/2]_+[s-1-u]_+[1][u]}
{[s+1/2]_+^2[1+u][u+3/2]},\\
&&R_{00}^{00}(u,s)=
\frac{[3+u][1][3/2-u]}{
[3][1+u][3/2+u]}+H_s
\frac{[1][u]}{[3][1+u]}.
\end{eqnarray*}
Here we have set
\begin{eqnarray}
G_s^\pm=-\frac{[2s\pm 2][s]_+}{[2s][s\pm 1]_+},~~~
H_s=G_s^+\frac{[s-5/2]_+}{[s+1/2]_+}+G_s^-
\frac{[s+5/2]_+}{[s-1/2]_+}.
\end{eqnarray}
The function $\rho^+(u)$
is given by
\begin{eqnarray}
{\rho}^+(u)&=&-qz^{\frac{1}{r}}\frac{
\{pq^2z\}
\{pq^3z\}
\{pq^3z\}
\{pq^4z\}
\{1/z\}
\{q/z\}
\{q^5/z\}
\{q^6/z\}
}{
\{pz\}
\{pqz\}
\{pq^5z\}
\{pq^6z\}
\{q^2/z\}
\{q^3/z\}
\{q^3/z\}
\{q^4/z\}
},\label{def:rhop}
\end{eqnarray}
where $z=q^{2u}$ and 
\begin{eqnarray}
\{z\}=(z;p,q^6)_\infty.
\end{eqnarray}
The $R$-matrix  $R^{+*}(u,s)=
(\pi_{V,z_1}\otimes \pi_{V,z_2})\cR(r^*,s)$
 is obtained from $R^{+}(u,s)$ by the 
replacements $r \to r^*$.
Hence, under the parametrization \eqref{sj},  
the dynamical $RLL$-relation takes the form 
\bea
&&R^{+(12)}(u,s+\al_1^\rv)
L^{+(1)}(u_1,s)L^{+(2)}(u_2,s+ \al_1^{\rv(1)})
=
L^{+(2)}(u_2,s)L^{+(1)}(u_1,s+\al_1^{\rv(2)})
R^{+*(12)}(u,s).\nn\\
\lb{DRLL2}
\ena

\subsection{Intertwining operators}

Let $\F, \F'$ be the highest weight $U_q$-modules.
We denote the type-I and type II intertwining operators of $U_q$-modules 
by $\Phi(z)$
and 
$\Psi^*(z)$, respectively. 
\begin{eqnarray}
\Phi(z):~{\cal F} \longrightarrow {\cal F}'\otimes W_z,~~\qquad
\Psi^*(z):~W_z \otimes {\cal F} \longrightarrow
{\cal F}'.
\end{eqnarray}
Twisting these operators by $F(r^*,s)$,
we obtain the corresponding intertwining operators ${\Phi}(v,s)$ 
and ${\Psi}^*(u,s)$ of $\Bqla$-modules.
\begin{eqnarray}
&&{\Phi}_W(u,s)=(\id \otimes \pi_{W,z})F(r^*,s) 
\Phi(z),\lb{intI}\\
&&{\Psi}_W^*(u,s)=
\Psi^*(z) (\pi_{W,z} \otimes \id)F(r^*,s)^{-1}.
\lb{intII}
\end{eqnarray}
From the intertwining relation satisfied by 
$\Phi(z)$ and $\Psi^*(z)$, 
one can derive  
the following dynamical intertwining relation for the  new intertwiners 
\cite{JKOS1}.
\begin{eqnarray}
&&{\Phi}_W^{(3)}(u_2+\frac{c}{2},s)
L_V^{+(1)}(u_1,s)=R^{+(13)}_{VW}(u,s+\al_1^\rv)
L_V^{+(1)}(u_1,s){\Phi}_W^{(3)}(u_2+\frac{c}{2},s+\al_1^{\rv(1)}),
\nonumber\\
\lb{dintrelI}\\
&&L_V^{+(1)}(u_1,s){\Psi}_W^{*(2)}(z_2,s+\al_1^{\rv(1)})=
{\Psi}_W^{*(2)}(z_2,s)L_V^{+(1)}
(u_1,s+\al_1^{\rv(2)})R_{VW}^{+*(12)}(u_1-u_2,s).
\nonumber
\\
\lb{dintrelII}
\end{eqnarray}
Note that \eqref{dintrelI} and \eqref{dintrelII} are the relations for
the operators $V_{z_1}\otimes {\cal F} \to 
V_{z_1}\otimes {\cal F} \otimes W_{z_2}$ and 
$V_{z_1}\otimes W_{z_2}\otimes {\cal F} \to 
V_{z_1}\otimes {\cal F}$, respectively. 


\section{Elliptic Algebra $U_{q,p}(A_2^{(2)})$}
In this section,
we give a definition of the elliptic algebra 
$U_{q,p}(A_2^{(2)})$.
We follows the idea given in
\cite{JKOS2, KoKo},
where the elliptic algebras $U_{q,p}(\g)$
 with non-twisted affine Lie algebra $\g$
are discussed.
We first introduce the  currents
$e(z,p),f(z,p)$ and $\psi^\pm(z,p)$
of the quantum group $U_q(A_2^{(2)})$,
by modifying the Drinfeld currents of $U_q(A_2^{(2)})$.
We then introduce the new current $k(z)$ in $U_q(A_2^{(2)})$
which is a more basic object than the currents $\psi^\pm(z,p)$.
Finally modifying them by taking a tensor product
with Heisenberg algebra, we introduce the elliptic currents
$E(u), F(u), H^\pm(u)$ and $K(u)$ forming the elliptic algebra
$U_{q,p}(A_2^{(2)})$. The current $K(u)$ plays an essential role 
in the $RLL$-formulation of $U_{q,p}(A_2^{(2)})$.
Hereafter we set $h=\al_1^\rv$.

\subsection{Drinfeld Currents of $U_q(A_2^{(2)})$}

Let us recall the Drinfeld currents
of $U_q(A_2^{(2)})$.
Let $0<q<1$.
We use the standard symbol of $q$-integer 
\begin{eqnarray}
[n]_q=\frac{q^n-q^{-n}}{q-q^{-1}}.
\end{eqnarray}

\begin{df}{\bf (Drinfeld currents)}~~
Let $x_m^\pm$ $(m \in {\mathbb{Z}})$,
$a_m (m \in {\mathbb{Z}}_{\neq 0})$
$q^c, q^h, q^d$
denote the Drinfeld generator of $U_q(A_2^{(2)})$.
In terms of the generating functions
\begin{eqnarray}
&&x^\pm(z)=\sum_{m \in {\Z}}x_m^\pm z^{-m},\\
&&\psi(q^{c/2}z)=q^{h/2}
\exp\left((q-q^{-1})\sum_{m>0}a_m z^{-m}\right)
,\\
&&\varphi(q^{-c/2}z)=q^{-h/2}
\exp\left(-(q-q^{-1})\sum_{m>0}a_{-m} z^{m}\right).
\end{eqnarray}
the defining relations of $U_q(A_2^{(2)})$
are given by
\begin{eqnarray}
&&q^c : {\rm central},~~~
q^d a_m q^{-d}=q^m a_m,~~
q^d q_m^\pm q^{-d}=q^m x_m^\pm,\\
&&q^{h}x^\pm(z)q^{-h}=
q^{\pm 2} x^\pm(z),~~q^dq^{h}=q^{h}q^d.
\end{eqnarray}
\begin{eqnarray}
&&~[a_m,a_n]=\delta_{m+n,0}
\frac{1}{m}([2m]_q-[m]_q)q^{-c|m|}[cm]_q,\\
&&~[a_m,x^+(z)]=
\frac{1}{m}([2m]_q-[m]_q)q^{-c|m|}z^m x^+(z),\\
&&~[a_m,x^-(z)]=
-\frac{1}{m}([2m]_q-[m]_q)z^m x^-(z),\\
&&(z_1-q^{\pm 2}z_2)
(z_1-q^{\mp 1}z_2)x^\pm(z_1)x^\pm(z_2)=
-(q^{\pm 2}z_1-z_2)
(q^{\mp 1}z_1-z_2)x^\pm(z_2)x^\pm(z_1),\nn\\
&&\\
&&~[x^+(z_1),x^-(z_2)]=\frac{1}{q-q^{-1}}
\left(\psi(q^{c/2}z_2)\delta(q^{-c}z_1/z_2)
-\varphi(q^{-c/2}z_2)\delta(q^cz_1/z_2)\right),\\
&&~\sum_{\sigma \in S_3}
\left(q^{\pm 3/2}z_{\sigma(1)}-
(q^{1/2}+q^{-1/2})z_{\sigma(2)}+q^{\pm 3/2}z_{\sigma(3)}
\right)x^\pm(z_{\sigma(1)})x^\pm(z_{\sigma(2)})
x^\pm(z_{\sigma(3)})=0.\nonumber\\
\end{eqnarray}
Here $\delta(z)$ denotes the delta function
$\delta(z)=\sum_{m \in {\mathbb{Z}}}z^m$.
We call the generators $h,a_m,x_m^\pm,c,d$
the Drinfeld generators of $U_q(A_2^{(2)})$
and the generating functions $x^\pm(z),\psi(z)$ and
$\varphi(z)$ the Drinfeld currents.
\end{df}
 

\subsection{Elliptic currents of $U_q(A_2^{(2)})$}

We next consider an elliptic modification of
the Drinfeld currents $x^\pm(z),\psi(z)$ and $\varphi(z)$.
Let us introduce the two auxiliary currents $u^\pm(z,p)$ by
\begin{eqnarray}
&&u^+(z,p)=\exp\left(
\sum_{m>0}\frac{a_{-m}}{[r^*m]_q}q^{rm}z^m
\right),\\
&&u^-(z,p)=\exp\left(
-\sum_{m>0}\frac{a_m}{[rm]_q}q^{rm}z^{-m}
\right).
\end{eqnarray}

\begin{prop}~~~The following commutation relations hold.
\begin{eqnarray}
&&u^+(z_1,p)u^-(z_2,p)\nn\\
&&\quad=u^-(z_2,p)u^+(z_1,p)
\frac{
(pq^{-c-2}z_1/z_2;p)_\infty
(p^*q^{c+2}z_1/z_2;p^*)_\infty
(pq^{-c+1}z_1/z_2;p)_\infty
(p^*q^{c-1}z_1/z_2;p^*)_\infty
}{
(pq^{-c+2}z_1/z_2;p)_\infty
(p^*q^{c-2}z_1/z_2;p^*)_\infty
(pq^{-c-1}z_1/z_2;p)_\infty
(p^*q^{c+1}z_1/z_2;p^*)_\infty
},\nn\\
&&\\
&&u^+(z_1,p)x^+(z_2)=
\frac{(p^*q^2z_1/z_2;p^*)_\infty 
(p^*q^{-1}z_1/z_2;p^*)_\infty
}{(p^*q^{-2}z_1/z_2;p^*)_\infty 
(p^*qz_1/z_2;p^*)_\infty}
x^+(z_2)u^+(z_1,p),\\
&&u^+(z_1,p)x^-(z_2)=
\frac{(p^*q^{c-2}z_1/z_2;p^*)_\infty
(p^*q^{c+1}z_1/z_2;p^*)_\infty
}
{
(p^*q^{c+2}z_1/z_2;p^*)_\infty
(p^*q^{c-1}z_1/z_2;p^*)_\infty
}x^-(z_2)u^+(z_1,p),\\
&&u^-(z_1,p)x^+(z_2)=
\frac{(pq^{-c-2}z_2/z_1;p)_\infty
(pq^{-c+1}z_2/z_1;p)_\infty
}{
(pq^{-c+2}z_2/z_1;p)_\infty
(pq^{-c-1}z_2/z_1;p)_\infty
}
x^+(z_2)u^-(z_1,p)
,\\
&&
u^-(z_1,p)x^-(z_2)=
\frac{(pq^2z_2/z_1;p)_\infty
(pq^{-1}z_2/z_1;p)_\infty
}{
(pq^{-2}z_2/z_1;p)_\infty
(pqz_2/z_1;p)_\infty
}x^-(z_2)u^-(z_1,p),\\
&&\psi(z_1,p)u^+(z_2,p)\nn\\
&&\quad=u^+(z_2,p)\psi(z_1,p)
\frac{(q^{r^*+2}z_2/z_1;p)_\infty
(q^{r^*-1}z_2/z_1;p)_\infty
(q^{r^*-2}z_2/z_1;p^*)_\infty
(q^{r^*+1}z_2/z_1;p^*)_\infty
}{(q^{r^*-2}z_2/z_1;p)_\infty
(q^{r^*+1}z_2/z_1;p)_\infty
(q^{r^*+2}z_2/z_1;p^*)_\infty
(q^{r^*-1}z_2/z_1;p^*)_\infty
},\\
&&\psi(z_1,p)u^-(z_2,p)\nn\\
&&\quad=u^-(z_2,p)\psi(z_1,p)
\frac{(q^{r-2}z_1/z_2;p)_\infty
(q^{r+1}z_1/z_2;p)_\infty
(q^{r+2}z_1/z_2;p^*)_\infty
(q^{r-1}z_1/z_2;p^*)_\infty
}{
(q^{r+2}z_1/z_2;p)_\infty
(q^{r-1}z_1/z_2;p)_\infty
(q^{r-2}z_1/z_2;p^*)_\infty
(q^{r+1}z_1/z_2;p^*)_\infty
}\nonumber
,\\
&&\psi(z_1,p)x^+(z_2)\nn\\
&&\quad=x^+(z_2)\psi(z_1,p)
\frac{(q^{r^*-2}z_2/z_1;p)_\infty
(q^{r^*+1}z_2/z_1;p)_\infty
(q^{r^*+2}z_1/z_2;p^*)_\infty
(q^{r^*-1}z_1/z_2;p^*)_\infty
}{(q^{r^*+2}z_2/z_1;p)_\infty
(q^{r^*-1}z_2/z_1;p)_\infty
(q^{r^*-2}z_1/z_2;p^*)_\infty
(q^{r^*+1}z_1/z_2;p^*)_\infty}
,\\
&&\psi(z_1,p)x^-(z_2)\nn\\
&&\quad=x^-(z_2)\psi(z_1,p)
\frac{(q^{r+2}z_2/z_1;p)_\infty 
(q^{r-1}z_2/z_1;p)_\infty (q^{r-2}z_1/z_2;p^*)_\infty 
(q^{r+1}z_1/z_2;p^*)_\infty}
{(q^{r-2}z_2/z_1;p)_\infty (q^{r+1}z_2/z_1;p)_\infty 
(q^{r+2}z_1/z_2;p^*)_\infty (q^{r-1}z_1/z_2;p^*)_\infty}.
\end{eqnarray}
\end{prop}

\begin{df}
We define ``dressed'' currents $e(z,p), f(z,p)$ and  $\psi^\pm(z,p)$ by
\begin{eqnarray}
&&e(z,p)=u^+(z,p)x^+(z),\\
&&f(z,p)=x^-(z)u^-(z,p),\\
&&\psi^+(z,p)=u^+(q^{c/2}z,p)\psi(z)u^-(q^{-c/2}z,p),\\
&&\psi^-(z,p)=u^+(q^{-c/2}z,p)\varphi(z)
u^-(q^{c/2}z,p).
\end{eqnarray}
\end{df}
If we introduce the auxiliary current $\psi(z,p)$ by
\begin{eqnarray}
\psi(z,p)=
\exp\left(\sum_{m>0}
\frac{x^{cm}}{[r^*m]_q}a_{-m}z^{m}
\right)
\exp\left(-\sum_{m>0}
\frac{1}{[rm]_q}a_{m}z^{-m}
\right),\lb{def:psi}
\end{eqnarray}
we have
\begin{eqnarray}
&&\psi^\pm(q^{\mp(r-c/2)}z)=
q^{\pm h/2}\psi(z,p).
\end{eqnarray}

\begin{prop}~~~
The currents 
$e(z,p), f(z,p)$ and $\psi(z,p)$
satisfy the following commutation relations.
\begin{eqnarray}
&&\psi(z_1,p)\psi(z_2,p)
=
\frac{\Theta_p(q^{-2}z_1/z_2)
\Theta_p(qz_1/z_2)}
{\Theta_p(q^2z_1/z_2)
\Theta_p(q^{-1}z_1/z_2)}
\frac{\Theta_{p^*}(q^2z_1/z_2)
\Theta_{p^*}(q^{-1}z_1/z_2)}
{\Theta_{p^*}(q^{-2}z_1/z_2)
\Theta_{p^*}(qz_1/z_2)}
\psi(z_2,p)\psi(z_1,p),\nn\\
&&\\
&&\psi(z_1,p)e(z_2,p)=
\frac{\Theta_{p^*}(q^{r^*+2}z_1/z_2)
\Theta_{p^*}(q^{r^*-1}z_1/z_2)
}{\Theta_{p^*}(q^{r^*-2}z_1/z_2)
\Theta_{p^*}(q^{r^*+1}z_1/z_2)}
e(z_2,p)\psi(z_1,p),\\
&&
\psi(z_1,p)f(z_2,p)=\frac{\Theta_{p}
(q^{r-2}z_1/z_2)
\Theta_{p}(q^{r+1}z_1/z_2)
}{\Theta_{p}(q^{r+2}z_1/z_2)
\Theta_{p}(q^{r-1}z_1/z_2)}
f(z_2,p)\psi(z_1,p),\\
&&
e(z_1,p)e(z_2,p)=(-1)
\frac{\Theta_{p^*}(q^{-2}z_2/z_1)
\Theta_{p^*}(q^{-1}z_1/z_2)}
{\Theta_{p^*}(q^{-2}z_1/z_2)
\Theta_{p^*}(q^{-1}z_2/z_1)}e(z_2,p)e(z_1,p),\\
&&
f(z_1,p)f(z_2,p)=(-1)
\frac{\Theta_p(q^{2}z_2/z_1)
\Theta_p(qz_1/z_2)}
{\Theta_p(q^{2}z_1/z_2)
\Theta_p(q z_2/z_1)}f(z_2,p)f(z_1,p),\\
&&~[e(z_1,p),f(z_2,p)]=
\frac{1}{q-q^{-1}}\left(
\psi^+(q^{c/2}z_2)\delta(q^{-c}z_1/z_2)
-\psi^-(q^{-c/2}z_2)\delta(q^cz_1/z_2)\right),
\end{eqnarray}
\begin{eqnarray}
&&\sum_{\sigma \in S_3}
\frac{(p^*q^2z_{\sigma(3)}/z_{\sigma(1)};p^*)_\infty
(p^*q^{-1}z_{\sigma(3)}/z_{\sigma(1)};p^*)_\infty
(p^*q^{-1}z_{\sigma(3)}/z_{\sigma(2)};p^*)_\infty
(p^*q^{-1}z_{\sigma(2)}/z_{\sigma(1)};p^*)_\infty
}{
(p^*q^{-2}z_{\sigma(3)}/z_{\sigma(1)};p^*)_\infty
(p^*q z_{\sigma(3)}/z_{\sigma(1)};p^*)_\infty
(p^*q z_{\sigma(3)}/z_{\sigma(2)};p^*)_\infty
(p^*q z_{\sigma(2)}/z_{\sigma(1)};p^*)_\infty
}
\nonumber\\
&&\quad\times
\left(
z_{\sigma(1)}\frac{
(q^2z_{\sigma(2)}/z_{\sigma(1)};p^*)_\infty
(p^*q^2 z_{\sigma(3)}/z_{\sigma(2)};p^*)_\infty
}{
(p^*q^{-2} z_{\sigma(2)}/z_{\sigma(1)};p^*)_\infty
(p^*q^{-2} z_{\sigma(3)}/z_{\sigma(2)};p^*)_\infty}\right.
\nonumber\\
&&\quad -\left.
qz_{\sigma(2)}\frac{
(p^*q^2 z_{\sigma(2)}/z_{\sigma(1)};p^*)_\infty
(p^*q^2 z_{\sigma(3)}/z_{\sigma(2)};p^*)_\infty}{
(p^*q^{-2}z_{\sigma(2)}/z_{\sigma(1)};p^*)_\infty
(p^*q^{-2}z_{\sigma(3)}/z_{\sigma(2)};p^*)_\infty}
\right)
e(z_{\sigma(1)},p)e(z_{\sigma(2)},p)
e(z_{\sigma(3)},p)=0,\nn\\
&&\\
&&\sum_{\sigma \in S_3}
\frac{(pqz_{\sigma(2)}/z_{\sigma(1)};p)_\infty
(pq^{-2}z_{\sigma(3)}/z_{\sigma(1)};p)_\infty
(pqz_{\sigma(3)}/z_{\sigma(2)};p)_\infty
(pqz_{\sigma(3)}/z_{\sigma(2)};p)_\infty
}{
(pq^{-1}z_{\sigma(2)}/z_{\sigma(1)};p)_\infty
(pq^2 z_{\sigma(3)}/z_{\sigma(1)};p)_\infty
(pq^{-1} z_{\sigma(3)}/z_{\sigma(1)};p)_\infty
(pq^{-1} z_{\sigma(3)}/z_{\sigma(2)};p)_\infty
}
\nonumber\\
&&\quad\times
\left(
z_{\sigma(1)}\frac{
(q^{-2}z_{\sigma(2)}/z_{\sigma(1)};p)_\infty
(pq^{-2} z_{\sigma(3)}/z_{\sigma(2)};p)_\infty
}{
(pq^{2} z_{\sigma(2)}/z_{\sigma(1)};p)_\infty
(pq^{2} z_{\sigma(3)}/z_{\sigma(2)};p)_\infty}\right.
\nonumber\\
&&\quad-\left.
q^{-1}z_{\sigma(2)}\frac{
(pq^{-2} z_{\sigma(2)}/z_{\sigma(1)};p)_\infty
(q^{-2} z_{\sigma(3)}/z_{\sigma(2)};p)_\infty}{
(pq^{2}z_{\sigma(2)}/z_{\sigma(1)};p)_\infty
(pq^{2}z_{\sigma(3)}/z_{\sigma(2)};p)_\infty}
\right)
f(z_{\sigma(1)},p)f(z_{\sigma(2)},p)
f(z_{\sigma(3)},p)=0.\nn\\
&&
\end{eqnarray}
\end{prop}

\subsection{Basic current $k(z)$}

The current $\psi(z,p)$ \eqref{def:psi} can be expressed by 
more basic current $k(z)$. 
Let us introduce new generator of the bosons,
$\alpha_m, \beta_m$.
\begin{eqnarray}
&&\alpha_m=\left\{\begin{array}{cc}
a_m,&  m>0,\\
\frac{[r m]_q}{[r^*m]_q} q^{c|m|} a_m,& m<0,
\end{array}\right.~~~~~\beta_m=\alpha_m 
\frac{[r^*m]_q}{[rm]_q}.
\end{eqnarray}
They satisfy the following commutation relations.
\begin{eqnarray}
&&~[\alpha_m,\alpha_n]=\delta_{m+n,0}
\frac{[2m]_q-[m]_q}{m}
\frac{[cm]_q [rm]_q}{[r^*m]_q},\\
&&~
[\beta_m,\beta_n]=\delta_{m+n,0}
\frac{[2m]_q-[m]_q}{m}
\frac{[cm]_q [r^*m]_q}{[r m]_q},
\end{eqnarray}
The current $\psi(z,p)$ is expressed by
\begin{eqnarray}
\psi(z,p)=:\exp\left(-\sum_{m \neq 0}
\frac{\alpha_m}{[rm]_q}z^{-m}\right):=
:\exp\left(-\sum_{m \neq 0}
\frac{\beta_m}{[r^*m]_q}z^{-m}\right):.
\end{eqnarray}
The colons $:~:$ denote the standard normal ordering.

\begin{df}~~{\bf (Basic Current)}~~
We define the current $k(z)$ by
\begin{eqnarray}
k(z,p)=:\exp\left(-\sum_{m \neq 0}\frac{[m]_q}{[rm]_q
([2m]_q-[m]_q)}
\al_m z^{-m}\right):.
\end{eqnarray}
\end{df}
The current $\psi(z,p)$ is expressed by $k(z)$ as follows.
\begin{eqnarray}
\psi(z,p)=:k(q^{-1}z,p)k(z,p)^{-1}k(qz,p):.
\end{eqnarray}
By a straightforward calculation,
we have the following commutation relations.
\begin{prop}
\begin{eqnarray}
&&k(z_1,p)u^+(z_2,p)=
\frac{(q^{r^*+1}z_2/z_1;p)_\infty
(q^{r^*-1}z_2/z_1;p^*)_\infty
}{(q^{r^*-1}z_2/z_1;p)_\infty
(q^{r^*+1}z_2/z_1;p^*)_\infty}
u^+(z_2,p)k(z_1,p),\\
&&k(z_1,p)u^-(z_2,p)=\frac{(q^{r-1}z_1/z_2;p)_\infty 
(q^{r+1}z_1/z_2;p^*)_\infty}
{(q^{r+1}z_1/z_2;p)_\infty
(q^{r-1}z_1/z_2;p^*)_\infty}
u^-(z_2,p)k(z_1,p)
,\\
&&k(z_1,p)x^+(z_2)=
\frac{(q^{r^*+1}z_1/z_2;p^*)_\infty 
(q^{r^*-1}z_2/z_1;p)_\infty}
{(q^{r^*-1}z_1/z_2;p^*)_\infty 
(q^{r^*+1}z_2/z_1;p)_\infty}
x^+(z_2)k(z_1,p)
,\\
&&k(z_1,p)x^-(z_2)=
\frac{(q^{r-1}z_1/z_2;p^*)_\infty
(q^{r+1}z_2/z_1;p)_\infty
}{
(q^{r+1}z_1/z_2;p^*)_\infty
(q^{r-1}z_2/z_1;p)_\infty
}
x^-(z_2)k(z_1,p).
\end{eqnarray}
\end{prop}

\begin{prop}~~The currents $e(z,p), f(z,p)$ and
$k(z,p)$ satisfy the following commutation relations.
\begin{eqnarray}
k(z_1,p)k(z_2,p)&=&
z^{-1/r^*+1/r}
\rho(z_1/z_2)
k(z_2,p)k(z_1,p),\\
k(z_1,p)e(z_2,p)&=&
\frac{\Theta_{p^*}(q^{r^*+1}z_1/z_2)}
{\Theta_{p^*}(q^{r^*-1}z_1/z_2)}e(z_2,p)k(z_1),\\
k(z_1,p)f(z_2,p)&=&
\frac{\Theta_{p}(q^{r-1}z_1/z_2)}
{\Theta_{p}(q^{r+1}z_1/z_2)}
f(z_2,p)k(z_1,p).
\end{eqnarray}
Here we have set
\begin{eqnarray}
{\rho}(z)&=&\frac{{\rho}^{+*}(z)}{
{\rho}^{+}(z)},\label{def:rho}
\end{eqnarray}
where $\rho^+(z)$ is given in (\ref{def:rhop}) and
$\rho^{+*}(z)=\rho^+(z)|_{r\to r^*}$.
\end{prop}

\subsection{Elliptic Algebra $U_{q,p}(A_2^{(2)})$}

Now we give a definition of the elliptic algebra
$U_{q,p}(A_2^{(2)})$. For this purpose, we 
 introduce a Heisenberg algebra $\mathbb{C}\{{\cal H}\}$ 
generated by $P,Q$, and $\bar{\alpha}$.
\begin{eqnarray}
&&~[P,Q]=1,~~~[Q,\bar{\alpha}]=\pi i,~~~[P,\bar{\alpha}]=0,\lb{def:h1}\\
&&~[P,P]=[Q,Q]=[\bar{\alpha},\bar{\alpha}]=0.\lb{def:h2}
\end{eqnarray}

\begin{df}~~{\bf (Elliptic Currents)}~~
We define the elliptic (total) currents $E(z), F(z)$ and $K(z)$ by
\begin{eqnarray}
E(z)&=&e(z)e^{\bar{\alpha}}e^{-Q}z^{-P/r^*},
\label{real:EC1}\\
F(z)&=&f(z)e^{-\bar{\alpha}}z^{P/r+h/2r},
\label{real:EC2}\\
K(z)&=&k(z)e^{-Q}z^{(1/r-1/r^*)P+h/2r}.
\label{real:EC3}
\end{eqnarray}
\end{df}
Let us introduce the auxiliary currents
$H^\pm(z)$ by
\begin{eqnarray}
&&H^\pm(z)=H(q^{\pm(r-c/2)}z),\label{def:EA7}\\
&&H(z)=\psi(z)e^{-Q}z^{(1/r-1/r^*)P+h/2r}=\kappa
K(qz)K(z)^{-1}K(q^{-1}z),
\end{eqnarray}
where 
\begin{eqnarray}
\kappa=\frac{\{pq^8\}
\{pq^5\}
\{pq^3\}
\{pq^4\}^2
\{p\}
\{p^*q^7\}^*
\{p^*q\}^*
\{p^*q^2\}^{*2}
\{p^*q^6\}^{*2}
}{
\{pq^7\}
\{pq\}
\{pq^2\}^2
\{pq^6\}^2
\{p^*\}^*
\{p^*q^8\}^*
\{p^*q^5\}^*
\{p^*q^3\}^{*}
\{p^*q^4\}^{*2}
}.\label{def:kappa}
\end{eqnarray}
From the commutation relations of the currents
$e(z,p), f(z,p)$ and $k(z,p)$, we can verify the following
relations.
\begin{thm}~~~The elliptic currents
$E(z), F(z)$ and $K(z)$ satisfy the following commutation
relations.
\begin{eqnarray}
K(z_1)K(z_2)&=&
\rho(z_1/z_2)K(z_2)K(z_1),\label{def:EA1}\\
K(z_1)E(z_2)&=&-\frac{[u_1-u_2+\frac{r^*+1}{2}]^*}
{[u_1-u_2+\frac{r^*-1}{2}]^*}E(z_2)K(z_1),
\label{def:EA2}
\\
K(z_1)F(z_2)&=&-\frac{[u_1-u_2+\frac{r-1}{2}]}
{[u_1-u_2+\frac{r+1}{2}]}F(z_2)K(z_1),
\label{def:EA3}\\
E(z_1)E(z_2)&=&-\frac{[u_1-u_2+1]^*
[u_1-u_2-\frac{1}{2}]^*}{
[u_1-u_2-1]^*
[u_1-u_2+\frac{1}{2}]^*
}E(z_2)E(z_1),\label{def:EA4}
\\
F(z_1)F(z_2)&=&-\frac{[u_1-u_2-1][u_1-u_2+\frac{1}{2}]}
{[u_1-u_2+1][u_1-u_2-\frac{1}{2}]}F(z_2)F(z_1),
\label{def:EA5}
\end{eqnarray}
\begin{eqnarray}
~[E(z_1),F(z_2)]&=&\frac{1}{q-q^{-1}}
\left(H^+(q^{c/2}z_2)\delta(q^{-c}z_1/z_2)-
H^-(q^{-c/2}z_2)\delta(q^cz_1/z_2)\right).
\label{def:EA6}
\end{eqnarray}
Here 
${\rho}(z)$ is given in (\ref{def:rho}).
The elliptic currents $E(z)$ and $F(z)$
satisfy the following Serre relations.
\begin{eqnarray}
&&\sum_{\sigma \in S_3}
\frac{(p^*q^2z_{\sigma(3)}/z_{\sigma(1)};p^*)_\infty
(p^*q^{-1}z_{\sigma(3)}/z_{\sigma(1)};p^*)_\infty
(p^*q^{-1}z_{\sigma(3)}/z_{\sigma(2)};p^*)_\infty
(p^*q^{-1}z_{\sigma(2)}/z_{\sigma(1)};p^*)_\infty
}{
(p^*q^{-2}z_{\sigma(3)}/z_{\sigma(1)};p^*)_\infty
(p^*q z_{\sigma(3)}/z_{\sigma(1)};p^*)_\infty
(p^*q z_{\sigma(3)}/z_{\sigma(2)};p^*)_\infty
(p^*q z_{\sigma(2)}/z_{\sigma(1)};p^*)_\infty
}
\nonumber\\
&&\quad\times
z_{\sigma(1)}^{-\frac{1}{2r^*}}
z_{\sigma(2)}^{-\frac{1}{r^*}}
\left(
z_{\sigma(1)}\frac{
(q^2z_{\sigma(2)}/z_{\sigma(1)};p^*)_\infty
(p^*q^2 z_{\sigma(3)}/z_{\sigma(2)};p^*)_\infty
}{
(p^*q^{-2} z_{\sigma(2)}/z_{\sigma(1)};p^*)_\infty
(p^*q^{-2} z_{\sigma(3)}/z_{\sigma(2)};p^*)_\infty}\right.
\nonumber\\
&&\quad-\left.
qz_{\sigma(2)}\frac{
(p^*q^2 z_{\sigma(2)}/z_{\sigma(1)};p^*)_\infty
(p^*q^2 z_{\sigma(3)}/z_{\sigma(2)};p^*)_\infty}{
(p^*q^{-2}z_{\sigma(2)}/z_{\sigma(1)};p^*)_\infty
(p^*q^{-2}z_{\sigma(3)}/z_{\sigma(2)};p^*)_\infty}
\right)
E(z_{\sigma(1)})E(z_{\sigma(2)})
E(z_{\sigma(3)})=0,
\label{def:EA8}
\end{eqnarray}
and
\begin{eqnarray}
&&\sum_{\sigma \in S_3}
\frac{(pqz_{\sigma(2)}/z_{\sigma(1)};p)_\infty
(pq^{-2}z_{\sigma(3)}/z_{\sigma(1)};p)_\infty
(pqz_{\sigma(3)}/z_{\sigma(2)};p)_\infty
(pqz_{\sigma(3)}/z_{\sigma(2)};p)_\infty
}{
(pq^{-1}z_{\sigma(2)}/z_{\sigma(1)};p)_\infty
(pq^2 z_{\sigma(3)}/z_{\sigma(1)};p)_\infty
(pq^{-1} z_{\sigma(3)}/z_{\sigma(1)};p)_\infty
(pq^{-1} z_{\sigma(3)}/z_{\sigma(2)};p)_\infty
}
\nonumber\\
&&\quad\times
z_{\sigma(1)}^{2/r} z_{\sigma(2)}^{1/r}
\left(
z_{\sigma(1)}\frac{
(q^{-2}z_{\sigma(2)}/z_{\sigma(1)};p)_\infty
(pq^{-2} z_{\sigma(3)}/z_{\sigma(2)};p)_\infty
}{
(pq^{2} z_{\sigma(2)}/z_{\sigma(1)};p)_\infty
(pq^{2} z_{\sigma(3)}/z_{\sigma(2)};p)_\infty}\right.
\nonumber\\
&&\quad-\left.
q^{-1}z_{\sigma(2)}\frac{
(pq^{-2} z_{\sigma(2)}/z_{\sigma(1)};p)_\infty
(q^{-2} z_{\sigma(3)}/z_{\sigma(2)};p)_\infty}{
(pq^{2}z_{\sigma(2)}/z_{\sigma(1)};p)_\infty
(pq^{2}z_{\sigma(3)}/z_{\sigma(2)};p)_\infty}
\right)
F(z_{\sigma(1)})F(z_{\sigma(2)})
F(z_{\sigma(3)})=0.\nonumber\\
\label{def:EA9}
\end{eqnarray}
\end{thm}

\begin{df}~~~
{\bf (Elliptic Algebra $U_{q,p}(A_2^{(2)})$)}~~
We define the elliptic algebra $U_{q,p}(A_2^{(2)})$
to be the associative algebra generated by the 
currents $E(z), F(z)$ and $K(z)$ satisfying the relations
(\ref{def:EA7})-(\ref{def:EA9}).
\end{df}

\begin{cor}~~~
The construction
of $E(z), F(z)$ and $K(z)$ given in
(\ref{real:EC1})-(\ref{real:EC3}) 
is a realization of the elliptic algebra
$U_{q,p}(A_2^{(2)})$ in terms of the Drinfeld
generator of the quantum group $U_q(A_2^{(2)})$ and the Heisenberg algebra 
$\C\{\H\}$.
\end{cor}

For later convenience, let us introduce auxiliary currents
$K_\epsilon(z), (\epsilon=0,\pm)$ by
\begin{eqnarray}
K_+(z)&=&K(q^{r-2}z)
=k(q^{r-2}z)e^{-Q}(q^{r-2}z)^{(1/r-1/r^*)P+h/2r},\\
K_0(z)&=&K(q^rz)^{-1}K(q^{r-1}z)
=k(q^{r}z)^{-1}k(q^{r-1}z)
q^{(1/r^*-1/r)P-h/2r},\\
K_-(z)&=&K(q^{r+1}z)^{-1}=
k(q^{r+1}z)^{-1}(q^{r+1}z)^{(1/r^*-1/r)P-h/2r}e^Q.
\end{eqnarray}
Then one can verify the following relations.
\begin{prop}
\begin{eqnarray}
&&K_+(z_1)E(z_2)=-\frac{[u_1-u_2+\frac{c-1}{2}]^*}
{[u_1-u_2+\frac{c-3}{2}]^*}E(z_2)K_+(z_1),
\label{rel:EA1}
\\
&&K_0(z_1)E(z_2)=\frac{[u_1-u_2+\frac{c}{2}]^*
[u_1-u_2+\frac{c-1}{2}]^*
}
{[u_1-u_2+\frac{c}{2}-1]^*
[u_1-u_2+\frac{c+1}{2}]^*
}E(z_2)K_0(z_1),
\label{rel:EA2}
\\
&&K_-(z_1)E(z_2)=-\frac{[u_1-u_2+\frac{c}{2}]^*}
{[u_1-u_2+\frac{c}{2}+1]^*}E(z_2)K_-(z_1),
\label{rel:EA3}
\\
&&K_+(z_1)F(z_2)=-\frac{[u_1-u_2-\frac{3}{2}]}
{[u_1-u_2-\frac{1}{2}]}F(z_2)K_+(z_1),
\label{rel:EA4}
\\
&&K_0(z_1)F(z_2)=\frac{[u_1-u_2-1]
[u_1-u_2+\frac{1}{2}]
}
{[u_1-u_2]
[u_1-u_2-\frac{1}{2}]
}F(z_2)K_0(z_1),
\label{rel:EA5}
\\
&&K_-(z_1)F(z_2)=-\frac{[u_1-u_2+1]}
{[u_1-u_2]}F(z_2)K_-(z_1),
\label{rel:EA6}
\\
&&H^\pm(q^{\pm c/2}z)=H(q^{\pm r}z)=\kappa K_-(z
)^{-1}K_0(z)=
\kappa' K_+(qz)
K_0(qz)^{-1}.
\label{rel:EA7}
\end{eqnarray}
Here 
$\mu$ is given in (\ref{def:kappa}) and $\kappa'$
is given by
\begin{eqnarray}
\kappa'=\frac{\{pq^{10}\}
\{pq^7\}
\{pq^5\}
\{pq^6\}^2
\{pq^2\}
\{p^*q^9\}^*
\{p^*q^3\}^*
\{p^*q^5\}^{*2}
\{p^*q^8\}^{*2}
}{
\{pq^9\}
\{pq^3\}
\{pq^5\}^2
\{pq^8\}^2
\{p^*q^2\}^*
\{p^*q^{10}\}^*
\{p^*q^7\}^*
\{p^*q^5\}^{*}
\{p^*q^6\}^{*2}
}.\label{def:kappa'}
\end{eqnarray}
\end{prop}

\section{Half Currents}

As a preparation for the $RLL$-formulation of the elliptic 
algebra $U_{p,q}(A_2^{(2)})$ in the next section, we here introduce 
the half currents, and investigate their commutation relations.

Let us first summarize the commutation relations between 
 the Heisenberg algebra $\C\{{\H}\}$ and the elliptic currents.
From (\ref{real:EC1})-(\ref{real:EC3}), 
we have the following relations.
\begin{prop}
\begin{eqnarray}
&&~[E(z),P]=E(z),~~[F(z),P+\frac{h}{2}]=F(z),\\
&&~[E(z),P+\frac{h}{2}]=0,~[F(z),P]=0,\\
&&~[K_+(z),P]=K_+(z)=[K_+(z),P+\frac{h}{2}],\\
&&~[K_0(z),P]=0=[K_0(z),P+\frac{h}{2}],\\
&&~[K_-(z),P]=-K_-(z)=[K_-(z),P+\frac{h}{2}].
\end{eqnarray}
\end{prop}

Now we define the half currents 
$E_{-,0}^+(u), 
E_{0,+}^+(u), 
E_{-,+}^+(u),
F_{0,-}^+(u),
F_{+,0}^+(u), 
F_{+,-}^+(u)$ and
$K_\epsilon^+(u) ,(\epsilon=0,\pm),$ by the following formulae.
\begin{df}~{\bf(Half Currents)}
\begin{eqnarray}
K_\epsilon^+(u)&=&K_\epsilon(z),~~~(\epsilon=0,\pm),
\label{def:HC1}\\
E_{-,0}^+(u)&=&a_{-0}^*
\oint_{C_{-0}^*} 
\frac{dz'}{2\pi iz'} E(z')\frac{[u-u'-P+\frac{c+1}{2}]_+^*[1]^*}
{[u-u'+\frac{c}{2}]^*
[P-\frac{1}{2}]_+^*},
\label{def:HC2}
\\
E_{0,+}^+(u)&=&
a_{0+}^*
\oint_{C_{0+}^*} 
\frac{dz'}{2\pi iz'} E(z')\frac{[u-u'-P+\frac{c}{2}]_+^*[1]^*}
{[u-u'+\frac{c-1}{2}]^*
[P-\frac{1}{2}]_+^*},
\label{def:HC3}\\
E_{-,+}^+(u)&=&a_{-+}^*\oint \oint_{C_{-+}^*} 
\frac{dz'}{z'}
\frac{dz''}{z''}
E(z')E(z'')
\frac{[1]^{*2}}{[P-\frac{1}{2}]_+^*[2P-2]^*
}
\nonumber\\ 
&\times&
\frac{[u-u'-2P+2+\frac{c}{2}]^*[u'-u''-P]_+^*}{
[u-u'+\frac{c}{2}]^*[u'-u''-\frac{1}{2}]^*},
\label{def:HC4}
\\
\nonumber\\
F_{0,-}^+(u)&=&a_{0-}
\oint_{C_{0-}} 
\frac{dz'}{2\pi iz'} F(z')\frac{[u-u'+P+\frac{h-1}{2}]_+[1]}
{[u-u']
[P+\frac{h-1}{2}]_+},
\label{def:HC5}
\\
F_{+,0}^+(z)&=&a_{+0}
\oint_{C_{+0}} 
\frac{dz'}{2\pi iz'} F(z')
\frac{[u-u'+P+\frac{h}{2}-1]_+ [1]}{[u-u'-\frac{1}{2}]
[P+\frac{h-1}{2}]_+},\label{def:HC6}
\\
F_{+,-}^+(u)&=&a_{+-}\oint \oint_{C_{+-}} 
\frac{dz'}{2\pi iz'} 
\frac{dz''}{2\pi iz''}
F(z')F(z'') \frac{[P+\frac{h}{2}-1]_+[1]^2}{
[P+\frac{h-3}{2}]_+[P+\frac{h}{2}-2]_+[2P+h-2]}\nonumber\\
&\times&
\frac{[u-u'+2P+h-3][u-u''+1][u'-u''+P+
\frac{h}{2}-1]_+}
{[u-u'][u-u''][u'-u''+\frac{1}{2}]}.
\label{def:HC7}
\end{eqnarray}
Here $C_{-0}^*$ is a simple closed contour that encircles
$pq^cz$ but not $q^cz$.
We abbreviate it as $C_{-0}^*: |p^*q^cz|<|z'|<|q^cz|$.
Similarly the others are given by
\begin{eqnarray}
C_{-0}^*  &:&  |p^*q^cz|<|z'|<|q^cz|,\\
C_{0+}^*  &:&  |p^*q^{c-1}z|<|z'|<|q^{c-1}z|,\\
C_{-+}^*  &:&  |p^*q^c|<|z'|<|q^cz|,~|p^*q^cz|<|z''|<|q^cz|,~
|p^*q z'|<|z''|<|q z'|,  \\
C_{0-}  &:&  |pz|<|z'|<|z|,\\
C_{0+}^*  &:&  |pq^{-1}z|<|z'|<|q^{-1}z|,\\
C_{+-}    &:&  |pz|<|z'|<|z|,~|pz|<|z''|<|z|,~
|pqz'|<|z''|<|qz'|.  
\end{eqnarray}
The constants $a_{-0}^*, a_{0+}^*, a_{-+}^*$,
$a_{0-}, a_{+0}, a_{+-}$ are chosen to satisfy
\begin{eqnarray}
\frac{\mu  a_{0-}a_{-0}^* [1]^*}{q-q^{-1}}=-1=
\frac{\mu'  a_{+0}a_{0+}^* [1]^*}{q-q^{-1}},
~~a_{+-}=a_{0-}a_{0-},~~
a_{-+}^*=a_{-0}^* a_{-0}^*.
\end{eqnarray}
\end{df}
We can verify the following commutation relations.
\begin{thm}
\label{thm:HC}
~~~The half currents satisfy the following relations.
\begin{eqnarray}
&&K_\pm^+(u_1)K_\pm^+(u_2)=
\rho(u)K_\pm^+(u_2)K_\pm^+(u_1),
\label{rel:HC1}
\\
&&K_0^+(u_1)K_0^+(u_2)=
\frac{\rho(u)
\rho(u)
}{
\rho(u+\frac{1}{2})\rho(
u-\frac{1}{2})}
K_0^+(u_2)K_0^+(u_1),
\label{rel:HC2}
\\
&&K_-^+(u_1)K_+^+(u_2)=K_+^+(u_2)K_-^+(u_1)
\rho(u)\frac{[u_1-u_2+1][u_1-u_2+\frac{3}{2}]
[u_1-u_2]^*[u_1-u_2+\frac{1}{2}]^*}{
[u_1-u_2][u_1-u_2+\frac{1}{2}]
[u_1-u_2+1]^*[u_1-u_2+\frac{3}{2}]^*},
\label{rel:HC3}\nn\\
&&\\
&&K_-^+(u_1)K_0^+(u_2)=
\rho(u)\frac{[u_1-u_2]^*[u_1-u_2+1]}{
[u_1-u_2+1]^*[u_1-u_2]}
K_0^+(u_2)K_-^+(u_1),
\label{rel:HC4}
\\
&&K_0^+(u_1)K_+^+(u_2)=
\rho(u)
\frac{[u_1-u_2]^*[u_1-u_2+1]}{
[u_1-u_2+1]^*[u_1-u_2]}
K_+^+(u_2)K_0^+(u_1),
\label{rel:HC5}\\
&&K_-^+(u_1)^{-1}E_{-,0}^+(u_2)K_-^+(u_1)
=-E_{-,0}^+(u_2)\frac{[u+1]^*}{[u]^*}
+E_{-,0}^+(u_1)\frac{[P+\frac{1}{2}+u]_+^*[1]^*}{
[P+1/2]_+^* [u]^*},
\label{rel:HC6}
\\
&&K_+^+(u_2)^{-1}E_{0,+}^+(u_1)K_+^+(u_2)
=-\frac{[u+1]^*}{[u]^*}E_{0,+}^+(u_1)
+\frac{[-P+\frac{1}{2}+u]_+^*[1]^*}{
[-P+\frac{1}{2}]_+^* [u]^*}
E_{0,+}^+(u_2),
\label{rel:HC7}
\\
&&K_-^+(u_1)^{-1}
E_{-,+}^+(u_2)K_-^+(u_1)\nn\\
&&\quad=E_{-,+}^+(u_2)\frac{[u+\frac{3}{2}]^*[u+1]^*}
{[u+\frac{1}{2}]^*[u]^*}
+K_-^+(u_1)^{-1}E_{-,0}^+(u_2)K_-^+(u_1)E_{-,0}^+(u_1)
\frac{[-P-1-u]^*_+[1]^*}{[P+\frac{1}{2}]^*_+[u+
\frac{1}{2}]^*}\nn
\\
&&\qquad-E_{-,+}(u_1)
\left(\frac{[-2P+1-u]^*[u+\frac{3}{2}]^*[1]^*}
{[-2P+1]^*[u+\frac{1}{2}]^*[u]^*}
+\frac{
[2P-2]^*[P]_+^*
[-2P-\frac{1}{2}-u]^*[1]^*}{
[2P]^*[P-1]_+^*
[-2P+1]^*[u+\frac{1}{2}]^*}\right),
\label{rel:HC8}\\
&&K_-^+(u_1)F_{0,-}^+(u_2)K_-(u_1)^{-1}=
-\frac{[u+1]}{[u]}F_{0,-}^+(u_2)
+\frac{[-P+\frac{-h+1}{2}+u]_+[1]}
{[-P+\frac{-h+1}{2}]_+[u]}F_{0,-}^+(u_2),
\label{rel:HC9}
\\
&&K_+^+(u_2)
F_{+,0}^+(u_1)K_+^+(u_2)^{-1}=
-F_{+,0}^+(u_1)
\frac{[u+1]}{[u]}
+F_{+,0}^+(u_2)
\frac{[P+\frac{h+1}{2}+u]_+[1]}
{[P+\frac{h+1}{2}]_+[u]},
\label{rel:HC10}
\\
&&K_-^+(u_1)F_{+,-}^+(u_2)K_-^+(u_1)^{-1}\nn\\
&&\quad =
\frac{[u+\frac{3}{2}][u+1]}{[
u+\frac{1}{2}
][u]}F_{+,-}^+(u_2)
-\frac{[P+\frac{h}{2}-1-u]_+[1]}
{[-P+\frac{-h+1}{2}][u+\frac{1}{2}]}
F_{0,-}^+(u_1)K_-^+(u_1)F_{0,-}^+(u_2)
K_-^+(u_1)^{-1}
\nn\\
&&\qquad-\left(\frac{[u+\frac{3}{2}]
[2P+h+1-u][1]}
{[u+\frac{1}{2}][u][2P+h+1]}+
\frac{[2P+h-\frac{1}{2}-u]
[1][2P+h+2][P+\frac{h}{2}]_+}
{[2P+h+1][2P+h][P+\frac{h}{2}+1]_+[u+\frac{1}{2}]}
\right)F_{+,-}^+(u_1),\nn\\
&&\label{rel:HC11}\\
&&[E_{-,0}^+(u_1),F_{0,-}^+(u_2)]\nn\\
&&\quad=-K_0^+(u_2)
\frac{[-P-\frac{1}{2}+u]_+^*[1]^*}
{[-P-\frac{1}{2}]_+^*[u]^*}K_-^+(u_2)^{-1}+
K_-^+(u_1)^{-1}
\frac{[-P+\frac{-h+1}{2}+u]_+[1]}
{[-P+\frac{-h+1}{2}]_+[u]}K_0^+(u_1)
,\label{rel:HC12}\\
&&[E_{0,+}^+(u_1),F_{+,0}^+(u_2)]\nn\\
&&\quad=
-K_0^+(u_2)^{-1}\frac{[P+\frac{h+1}{2}+u]_+[1]}
{[P+\frac{h+1}{2}]_+[u+1]}K_+^+(u_2)
+K_+^+(u_1)\frac{[P-\frac{1}{2}+u]^*_+[1]^*}
{[P-\frac{1}{2}]^*_+[u+1]^*}
K_0^+(u_1)^{-1}.\label{rel:HC13}
\end{eqnarray}
where $u=u_1-u_2$.
\end{thm}
{\it Proof.}~~~
The relations (\ref{rel:HC1})-(\ref{rel:HC5})
are direct consequences of the commutation relation
of the elliptic current $K(z)$.
Let us consider the relations 
(\ref{rel:HC6})-(\ref{rel:HC13}).
These relations can be proved by
reducing them to identities of the theta functions.
We show the relation (\ref{rel:HC6}).
The relations (\ref{rel:HC7}), (\ref{rel:HC9}) and
(\ref{rel:HC10}) can be proved in the same way.
From the definition of the half current (\ref{def:HC2})
and the commutation relation of (\ref{rel:EA3}),
the LHS of \eqref{rel:HC6} yields 
\begin{eqnarray}
&&K_-^+(u_1)E_{-,0}^+(u_2)K_-^+(u_1)^{-1}\nonumber\\
&&\quad=-a_{-,0}^*
\oint_{C_{-0}^*}\frac{dz'}{2\pi iz'}
E(z')\frac{[u_1-u'+\frac{c}{2}+1]^*
[u_2-u'-P+\frac{c-1}{2}]_+^*[1]^*}
{[u_1-u'+\frac{c}{2}]^*[u_2-u'+\frac{c}{2}]^*
[P+\frac{1}{2}]_+^*}.
\end{eqnarray}
Then the equality is verified by the following identity of the theta functions.
\begin{eqnarray}
&-&\frac{[u_1-u'+\frac{c}{2}+1]^*
[u_2-u'-P+\frac{c-1}{2}]^*_+}
{[u_1-u'+\frac{c}{2}]^*
[u_2-u'+\frac{c}{2}]^*
[P+\frac{1}{2}]_+^*}\nn\\
&=&-\frac{[u_2-u'-P+\frac{c+1}{2}]_+^*
[u_1-u_2+1]^*}
{[u_2-u'+\frac{c}{2}]^*
[u_1-u_2]^*[P-\frac{1}{2}]_+^*}+
\frac{[u_1-u'-P+\frac{c+1}{2}]_+^*
[u_1-u_2+P+\frac{1}{2}]_+^*
[1]^*}
{[u_1-u'+\frac{c}{2}]^*
[u_1-u_2]^*[P-\frac{1}{2}]_+^*
[P+\frac{1}{2}]_+^*}.\nn
\end{eqnarray}

Next, we show the relation (\ref{rel:HC12}).
The relation (\ref{rel:HC13})
can be proved in the same way.
Integrating the delta function appearing in
(\ref{def:EA6}) and using (\ref{rel:EA7}),
we have
\begin{eqnarray}
&&(\mu a_{0-} a_{-0}^*)^{-1}(q-q^{-1})
[E_{-,0}^+(u_1),
F_{0,-}^+(u_2)]\nonumber\\
&=&
\oint_{C^+}\frac{dz'}{2\pi iz'}
K_-^+(u')^{-1}K_0^+(u')
\frac{[u_1-u'-P+\frac{1}{2}]_+^*
[1]^*
[u_2-u'+P+\frac{h-1}{2}]_+
[1]}
{[u_1-u']^*[P-\frac{1}{2}]_+^*
[u_2-u'][P+\frac{h-1}{2}]_+}
\\
&-&
\oint_{C^-}\frac{dz'}{2\pi iz'}
K_-^+(u'-r)^{-1}K_0^+(u'-r)
\frac{[u_1-u'-P+c+\frac{1}{2}]_+^*
[1]^*
[u_2-u'+P+\frac{h-1}{2}]_+
[1]}
{[u_1-u'+c]^*[P-\frac{1}{2}]_+^*
[u_2-u'][P+\frac{h-1}{2}]_+}\nonumber
\end{eqnarray}
Here the contours $C^\pm$ are now
\begin{eqnarray}
C^+:~|p^*z_1|,|pz_2|<|z'|<|z_1|, |z_2|,\\
C^-:~|pz_1|,|pz_2|<|z'|<|q^{2c}z_1|, |z_2|.
\end{eqnarray}
When we change the integration variable $z' \to pz'$
in the second term, 
the integrand becomes the same as the first term,
but 
the contour
$C^-$ is changed to $\tilde{C}^-$ given by 
\begin{eqnarray}
\tilde{C^-}:~|z_1|,|z_2|<|z'|<|p^{-1}q^{2c}z_1|, |p^{-1}z_2|.
\end{eqnarray}
Taking the residues at $z'=z_1, z_2$,
we get (\ref{rel:HC12}).

We give a proof of (\ref{rel:HC8}) in Appendix \ref{Proof}.
One can prove (\ref{rel:HC11}) in a similar way.
\begin{flushright}
{Q.E.D.}
\end{flushright}

\section{The $L$-operator of $U_{q,p}(A_2^{(2)})$
and Relation to ${\cal B}_{q,\lambda}(A_2^{(2)})$
}

In this section, we clarify the relation between
two elliptic algebras
$U_{q,p}(A_2^{(2)})$ and ${\cal B}_{q,\lambda}(A_2^{(2)})$.
For this purpose, we first construct a $L$-operator which gives 
the $RLL$-formulation of $U_{q,p}(A_2^{(2)})$. Then  
modifying $L$-operator by removing the Heisenberg generators $Q, \ba$, 
we derive the dynamical $RLL$-relation (\ref{dynRLL}) characterizing 
the elliptic quantum group ${\cal B}_{q,\lambda}(A_2^{(2)})$.

\subsection{$L$-operator of $U_{q,p}(A_2^{(2)})$}

\begin{df}
~~~By using the half currents, we define the $L$-operator
$\widehat{L}^+(u) \in {\rm End}(\C^3)\otimes U_{q,p}(A_2^{(2)})$
as follows.
\begin{eqnarray}
\widehat{L}^+(u)
&=&\left(
\begin{array}{ccc}
1&F_{+0}^+(u)&F_{+-}^+(u)\\
0&1&F_{0-}^+(u)\\
0&0&1
\end{array}
\right)
\left(
\begin{array}{ccc}
K_+^+(u)&0&0\\
0&K_0^+(u)&0\\
0&0&K_-^+(u)
\end{array}
\right)
\left(
\begin{array}{ccc}
1&0&0\\
E_{0+}^+(u)&1&0\\
E_{-+}^+(u)&E_{-0}^+(u)&1
\end{array}
\right).\nonumber\\
\label{def:L}
\end{eqnarray}
Here matrix elements are the half currents given in the previous
section.
\end{df}

By a direct comparison with the relations 
of the half currents 
appeared in Theorem \ref{thm:HC}, 
we get the
following commutation relations of the $L$-operator.

\noindent
\begin{thm}\lb{main}~~
The $L$-operator $\widehat{L}^+(u)$ satisfies
the following $RLL=LLR^*$ relation.
\begin{eqnarray}
&&R^{+(12)}(u_1-u_2,P+h/2)\hL^{+(1)}(u_1)\hL^{+(2)}(u_2)\nonumber\\
&&\qquad\qquad =\hL^{+(2)}(u_2)\hL^{+(1)}(u_1)R^{+*(12)}
(u_1-u_2,P-(h^{(1)}+h^{(2)})/2).\label{RLL}
\end{eqnarray}
\end{thm}
The above equation should be understood as
equation of 
the operators acting on the space $\C^3 \otimes \C^3
\otimes U_{q,p}(A_2^{(2)})$.
The operator $h$ in LHS acts on $U_{q,p}(A_2^{(2)})$,
whereas the operator $h^{(1)}+h^{(2)}$ in RHS
acts on $\C^3 \otimes \C^3$ as
$\left(\begin{array}{ccc}
2&0&0\\
0&0&0\\
0&0&-2
\end{array}\right)
\otimes 1+
1 \otimes
\left(\begin{array}{ccc}
2&0&0\\
0&0&0\\
0&0&-2
\end{array}\right)
$.

\subsection{
$U_{q,p}(A_2^{(2)})$
and 
${\cal B}_{q,\lambda}(A_2^{(2)})$
}

Based on the above theorem, 
we give a relation between $U_{q,p}(A_2^{(2)})$ and 
$\Bqla(A_2^{(2)})$. We argue that 
the $RLL$ relation \eqref{RLL} is equivalent to
the dynamical $RLL$ relation of $\Bqla(A_2^{(2)})$. 
Hence we can regard the elliptic currents 
in $U_{q,p}(A_2^{(2)})$ as 
an elliptic analogue 
of the Drinfeld currents in $U_q(A_2^{(2)})$  providing 
a new realization of the elliptic quantum group 
$\Bqla(A_2^{(2)})$.
In order to show this, 
we consider the realization of $U_{q,p}(A_2^{(2)})$
given in \eqref{real:EC1}-\eqref{real:EC3} and 
modify the half currents in such a way that they have no 
$Q, \bar{\alpha}$ dependence. 
\begin{eqnarray}
&&k_+(u,P)=K_+(u)e^Q,~~
k_0(u,P)=K_0(u)~~
k_-(u,P)=K_-(u)e^{-Q},\\
&&f_{+,0}(u,P)=e^{\bar{\alpha}}F_{+,-}(u),
~~
f_{0,-}(u,P)=e^{\bar{\alpha}}F_{0,-}(u),~~
f_{+,-}(u,P)=e^{\bar{\alpha}}F_{+,-}(u)e^{\bar{\alpha}},
\\
&&e_{0,+}(u,P)=E_{0,+}(u)e^Qe^{-\bar{\alpha}},~~
e_{-,0}(u,P)=e^Qe^{-\bar{\alpha}}E_{-,0}(u),~~
e_{-,+}(u,P)=
e^Qe^{-\bar{\alpha}}E_{-,+}(u)e^Qe^{-\bar{\alpha}}.\nonumber\\
\end{eqnarray}
We regard them as the currents in $U_q(A_2^{(2)})$ with
parameters $P$ and $r$. Then we define a dynamical $L$-operator
$\widehat{L}^+(u,P)$ by
\begin{eqnarray}
\widehat{L}^+(u,P)
&=&\left(
\begin{array}{ccc}
1&f_{+0}^+(u,P)&f_{+-}^+(u,P)\\
0&1&f_{0-}^+(u,P)\\
0&0&1
\end{array}
\right)\nonumber\\
&\times&
\left(
\begin{array}{ccc}
k_+^+(u,P)&0&0\\
0&k_0^+(u,P)&0\\
0&0&k_-^+(u,P)
\end{array}
\right)
\left(
\begin{array}{ccc}
1&0&0\\
e_{0+}^+(u,P)&1&0\\
e_{-+}^+(u,P)&e_{-0}^+(u,P)&1
\end{array}
\right).\nonumber\\
\label{def:modL}
\end{eqnarray}
The two $L$-operators $\widehat{L}^+(u)$ and 
$\widehat{L}^+(u,P)$ are related by
\begin{eqnarray}
\widehat{L}^+(u,P)=\widehat{L}^+(u)\left(\begin{array}{ccc}
e^{Q}&&\\
&1&\\
&&e^{-Q}
\end{array}\right)=
\widehat{L}^+(u)e^{Qh^{(1)}/2}.
\end{eqnarray}
Here $h^{(1)}=\left(\begin{array}{ccc}
2&0&0\\
0&0&0\\
0&0&-2
\end{array}\right)\otimes 1$.
Substituting this into (\ref{RLL}) and moving the factor $e^{-Qh^{(j)}/2}\ 
(j=1,2) $
to the right end in the both sides,
we get the following statement.
\begin{cor}~~~The dynamical $L$-operator 
$L^{+}(u,P)$ satisfies the dynamical $RLL$ relation.
\begin{eqnarray}
&&R^{+(12)}(u_1-u_2,P+h/2)L^{+(1)}(u_1,P)L^{+(2)}(u_2,P+h^{(1)}/2)
\nonumber\\
&&\qquad\qquad=L^{+(2)}(u_2,P)L^{+(1)}(u_1,P+h^{(2)}/2)R^{+*(12)}(u_1-u_2,P).
\label{dynRLL}
\end{eqnarray}
\end{cor}

Comparing this with \eqref{DRLL2}, 
we identify 
our $L^+(u,P)$ with $L^+(u,s)$ in \eqref{DRLL2} and 
$s$ with $P$. 
We hence regard the elliptic currents 
$E(z), F(z)$ and 
$K(z)$ in $U_{q,p}(A_2^{(2)})$ 
as the Drinfeld currents of the elliptic 
quantum group $\Bqla(A_2^{(2)})$, although  
$U_{q,p}(A_2^{(2)})$ and $\Bqla(A_2^{(2)})$ are different by tensoring
 the Heisenberg algebra $\C\{{\H}\}$. 
More precisely,  $U_{q,p}(A_2^{(2)})$ is  
an extension of the algebra $\Bqla(A_2^{(2)})$ 
by tensoring the Heisenberg algebra $\C\{{\H}\}$; 
first tensoring the generators $e^{Q},e^{\bar{\alpha}}$, then  
regarding $s=P$ and imposing  
the commutation relations \eqref{def:h1} and \eqref{def:h2}. 
Naively we regards
$U_{q,p}(A_2^{(2)})$ as $\Bqla(A_2^{(2)})
\otimes_{{\mathbb{C}}\{P\}}
{\mathbb{C}}\{{\H}\}$.


\section{Vertex Operators}

Tensoring the Heisenberg algebra breaks down 
the  coalgebra structure of $\Bqla(A_2^{(2)})$. But 
we can define the $U_{q,p}(A_2^{(2)})$ counterparts of 
the intertwining operators of $\Bqla(A_2^{(2)})$. 
We call such operators the vertex operators of 
$U_{q,p}(A_2^{(2)})$.
In this section, we study
such vertex operators 
and compare them with those of
the dilute $A_L$ model obtained in \cite{HJKOS}.

\subsection{Intertwining Relations}

We here derive the $U_{q,p}(A_2^{(2)})$ 
counterparts of the dynamical 
intertwining relations \eqref{dintrelI}-\eqref{dintrelII}. 
In the next subsection, we use such relations to derive
a free field realization of the vertex operators. 

Let us first define an extension of the $U_q$ modules by
\be
&&\widehat{\cal F}=\bigoplus_{\mu \in {\mathbb{Z}}}
{\cal F}\otimes e^{\mu Q}.
\en 
Let ${\Phi}_W(u,P)$ and ${\Psi}^*_W(u,P)$ be the 
type I and type II intertwining operators of 
$\Bqla(A_2^{(2)})$ \eqref{dintrelI}-\eqref{dintrelII}.
We define type I and type II vertex operators 
$\widehat{\Phi}_W(u),\
 \widehat{\Psi}^*_W(u)$ of $U_{q,p}(A_2^{(2)})$
as the following extensions of the corresponding 
intertwining operators of $\Bqla(A_2^{(2)})$.
\begin{eqnarray}
&&\widehat{\Phi}_W(u)={\Phi}_W(u+c/2,P)\qquad :
\widehat{\cal F} \longrightarrow \widehat{\cal F}' \otimes W_z,\\
&&\widehat{\Psi}^*_W(u)={\Psi}^*_W(u,P)
e^{h^{(1)}Q/{2}}\qquad :
W_z \otimes \widehat{\cal F} \longrightarrow \widehat{\cal F}'.
\end{eqnarray}
From the commutation relation of $P$ and $Q$,
the new operators $\widehat{\Phi}_W(u)$
and $\widehat{\Psi}_W^*(u)$
satisfy the following "intertwining relations". 
\begin{eqnarray}
\widehat{\Phi}^{(3)}_W(u_2)
\widehat{L}^{+(1)}_V(u_1)&=&R^{+(1,3)}_{VW}(u_1-u_2,P+h/2)
\widehat{L}^{+(1)}_V(u_1)\widehat{\Phi}^{(3)}_W(u_2),
\label{def:Type-I}\\
\widehat{L}^{+(1)}_V(u_1)\widehat{\Psi}^{*(2)}_W(u_2)&=&
\widehat{\Psi}^{*(2)}_W(u_2)
\widehat{L}^{+(1)}_V(u_1)
R^{+*(1,2)}_{VW}(u_1-u_2,P-(h^{(1)}+h^{(2)})/2).
\label{def:Type-II}
\end{eqnarray}
Now let us restrict ourselves to the vector representation
$W=V\cong \C v_+ \oplus \C v_0 \oplus \C v_{-}$.
In this case, the $R$-matrix $R^+_{VV}(u,P)$
is given by
$R^+(u,P)$ in (\ref{def:R}), and the $L$-operator
$\widehat{L}^+_{V}(u,P)$ by $\widehat{L}^+(u,P)$ in
(\ref{def:modL}).
Let us set the components of the vertex operators 
$\Phi_j(u), \Psi_j^*(u), (j=\pm,0)$ by
\begin{eqnarray}
&&\widehat{\Phi}\left(u-\frac{1}{2}\right)=
\sum_{j=\pm ,0}
\Phi_j(u)\otimes v_j,\\
&&\widehat{\Psi}^*\left(u-\frac{c+1}{2}\right)
(v_j\otimes \cdot)=
\Psi_j^*(u),
\end{eqnarray}
and the matrix elements of
the $L$-operator $\widehat{L}^+(u)$ by
\begin{eqnarray}
\widehat{L}^+(u)v_{j}=
\sum_{m=0,\pm}
v_m L^+_{m,j}(u).
\end{eqnarray}

Let us investigate
the relations 
 (\ref{def:Type-I}) and (\ref{def:Type-II})
in detail. From the components
$[~(-,-), (j)], j=\pm,0$ of (\ref{def:Type-I}),
we have 
\begin{eqnarray}
\Phi_-\left(u_2+\frac{1}{2}\right)L_{-,j}^+(u_1)=
\rho^+(u_1-u_2)L_{-,j}^+(u_1)\Phi_-\left(u_2+\frac{1}{2}\right).
\label{rel:Type-I1}
\end{eqnarray}
Putting the definition $L_{-,j}^+(u)=K_-^+(u)E_{-,j}^+(u)$ into
the above, we have 
\begin{eqnarray}
&&\Phi_-\left(u_2+\frac{1}{2}\right)
K_-^+(u_1)=\rho^+(u_1-u_2)K_-^+(u_1)\Phi_-
\left(u_2+\frac{1}{2}\right),
\label{rel:Type-I2}\\
&&[\Phi_-(u_1),E_{-,0}^+(u_2)]=0,
\label{rel:Type-I3}\\
&&[\Phi_-(u_1),E_{-,+}^+(u_2)]=0.
\label{rel:Type-I4}
\end{eqnarray}
We have the sufficient condition of 
(\ref{rel:Type-I3}), (\ref{rel:Type-I4}).
\begin{eqnarray}
&&\Phi_-(u_1)E(u_2)=E(u_2)\Phi_-(u_1),
~~[\Phi_-(z_1),P]=0.
\label{rel:Type-I5}
\end{eqnarray}
From the component $~[(0,-), (-)]$ of
(\ref{def:Type-I}), we have
\begin{eqnarray}
\Phi_-\left(u_2+\frac{1}{2}\right)
F_{0,-}^+(u_1)K_-^+(u_1)
&=&\rho^+(u)\bR_{0-}^{0-}(u,P+h/2)
F_{0,-}^+(u_1)K_-^+(u_1)\Phi_-\left(u_2+\frac{1}{2}\right)
\nonumber\\
&+&\rho^+(u)\bR_{0-}^{-0}(u,P+h/2)K_-^+(u_1)\Phi_0
\left(u_2+\frac{1}{2}\right).\label{rel:Type-I6}
\end{eqnarray}
Let us assume 
the operator product
$K_-^+(u_1)
\Phi_-\left(u_2+\frac{1}{2}\right)$
has no pole at $u_1-u_2=-1-r$.
Later we will check that, for $c=1$,
this assumption is satisfied in a free field realization.
Then from (\ref{rel:Type-I2}),
we conclude the operator product
$\Phi_-\left(u_2+\frac{1}{2}\right)
K_-^+(u_1)$ has zero at $u_1-u_2=-1-r$.
Therefore setting $u_1-u_2=-1-r$ in
(\ref{rel:Type-I6}), we have
\begin{eqnarray}
0=\frac{[P+h/2+1/2]_+}{[P+h/2-1/2]_+}
F_{0,-}^+(u_1)K_-^+(u_1)\Phi_-\left(u_2+\frac{1}{2}\right)
+K_-^+(u_1)\Phi_0\left(u_2+\frac{1}{2}\right).
\end{eqnarray}
Then we have
\begin{eqnarray}
\Phi_0(u)&=&F_{0,-}^+\left(u-r-\frac{1}{2}\right)
\Phi_-(u).\label{rel:Type-I7}
\end{eqnarray}
Substituting (\ref{rel:Type-I1}) and (\ref{rel:Type-I7})
into (\ref{rel:Type-I6}),
we get
\begin{eqnarray}
\Phi_-(u_1)F(u_2)&=&-\frac{[u_1-u_2+1/2]}
{[u_1-u_2-1/2]}F(u_2)\Phi_-(u_1).
\label{rel:Type-I8}
\end{eqnarray}

Similarly, in order to investigate the structure of the component $\Phi_+(u)$,
we have, from the components of $~[(+,-), (-)]~$ of
(\ref{def:Type-I}), 
\begin{eqnarray}
&&\Phi_-
\left(u_2+\frac{1}{2}\right)
F_{+,-}^+(u_1)K_-^+(u_1)\nn\\
&&\qquad=
\rho^+(u)\bR_{+-}^{+-}(u,P+h/2)
F_{+,-}^+(u_1)K_-^+(u_1)\Phi_-
\left(u_2+\frac{1}{2}\right)\nn\\
&&\qquad\quad+\rho^+(u)
\bR_{+-}^{00}(u,P+h/2)
F_{0,-}^+(u_1)K_-^+(u_1)\Phi_0
\left(u_2+\frac{1}{2}\right)\nonumber\\
&&\qquad\quad+\rho^+(u)\bR_{+-}^{-+}(u,P+h/2)
K_-^+(u_1)\Phi_+
\left(u_2+\frac{1}{2}\right).\label{rel:Type-I9}
\end{eqnarray}
On the other hands, 
from the component 
$[(+,-), (-,-)]$ of $RLL$ relation (\ref{RLL}), 
we have
\begin{eqnarray}
K_-^+(u_2)F_{+,-}^+(u_1)K_-^+(u_2)^{-1}&=&
R_{+-}^{+-}(u)F_{+,-}^+(u_1)+
R_{+-}^{00}(u)F_{0,-}^+(u_1)K_-^+(u_1)
F_{0,-}^+(u_2)K_-^+(u_1)^{-1}\nonumber\\
&&\qquad +
R_{+-}^{-+}(u)K_-^+(u_1)F_{+,-}^+(u_2)K_-^+(u_1)^{-1}.
\end{eqnarray}
Putting the above into (
\ref{rel:Type-I9}), we get
\begin{eqnarray}
&&\Phi_-\left(u_2+\frac{1}{2}\right)
F_{+,-}^+(u_1)K_-^+(u_1)\nonumber\\
&&\qquad=\rho^+(u)\bR_{+-}^{-+}(u|P+h/2)
K_-(u_1)\Phi_+\left(u_2+\frac{1}{2}\right)\nn\\
&&\qquad\quad+
\rho^+(u)K_-^+(u_2-r)F_{+,-}^+(u_1)
K_-^+(u_2-r)^{-1}
K_-^+(u_1)\Phi_-
\left(u_2+\frac{1}{2}\right)\nonumber\\
&&\qquad\quad+\rho^+(u)
\bR_{+-}^{-+}(u|P+h/2)
K_-^+(u_1)F_{+,-}^+(u_2-r)\Phi_-
\left(u_2+\frac{1}{2}\right).
\label{rel:Type-I10}
\end{eqnarray}
Note that at the point $u_1-u_2=-1-r$,
$\rho^+(u)$ has a zero, but $\rho^+(u)\bR_{+-}^{-+}(u|P+h/2)$
have no zeros. In addition, under the same assumption given just 
below \eqref{rel:Type-I6},  the product 
$\Phi_-\left(u_2+\frac{1}{2}\right)K_-^+(u_1)$ vanishes at $u_1-u_2=-1-r$.
Setting $u_1-u_2=-1-r$ in
(\ref{rel:Type-I10}), we thus have the following formula for $\Phi_+(z)$.  
\begin{eqnarray}
\Phi_+(u)&=&-F_{+,-}^+
\left(u-r-\frac{1}{2}\right)\Phi_-(u).
\label{rel:Type-I11}
\end{eqnarray}
In the next section, we construct a free field realization
of the type-I vertex operators using the relations
(\ref{rel:Type-I2}),
(\ref{rel:Type-I5}),
(\ref{rel:Type-I7}),
(\ref{rel:Type-I8}) and
(\ref{rel:Type-I11})
for $c=1$.
We can check that the resultant vertex operators
satisfy the intertwining relation
(\ref{def:Type-I}).

Similarly, the sufficient conditions for
the type-II vertex operators are extracted from
(\ref{def:Type-II}) as follows.
\begin{eqnarray}
&&\Psi_-^*\left(u_2+\frac{1+c}{2}\right)
K_-(u_1)\rho^{+*}(u)
=K_-^+(u_1)\Psi_-^*\left(u_2+\frac{1+c}{2}\right),
\label{rel:Type-II1}\\
&&\Psi_-^*(u_1)F(u_2)
=F(u_2)\Psi_-^*(u_1),~~
[\Psi_-^*(u),P+h/2]=0,
\label{rel:Type-II2}
\\
&&\Psi_-^*(u_1)E(u_2)=-\frac{[u_1-u_2-\frac{1}{2}]^*}
{[u_1-u_2+\frac{1}{2}]^*}E(u_2)\Psi_-^*(u_1),
\label{rel:Type-II3}\\
&&\Psi_+^*(u)=-\Psi_-^*(u)E_{-,+}
\left(u-\frac{1+c}{2}-r^*
\right),
\label{rel:Type-II4}
\\
&&\Psi_0^*(u)=\Psi_-^*(u)
E_{-,0}
\left(u-\frac{1+c}{2}-r^*
\right).
\label{rel:Type-II5}
\end{eqnarray}


\subsection{Free Field Realizations}

Now we construct a free field realization of the vertex
operators fixing the representation level $c=1$.
For this purpose, we introduce the simple root
operator $\alpha$, defined by
\begin{eqnarray}
~[h,\alpha]=2,~[a_m,\alpha]=0,~[P,\alpha]=0,~
[Q,\alpha]=0,~[\alpha, \bar{\alpha}]=0.
\end{eqnarray}
If we introduce $\hat{\alpha}$ by
\begin{eqnarray}
\hat{\alpha}=\alpha+\bar{\alpha}, 
\end{eqnarray}
we have
\begin{eqnarray}
~[h,\hat{\alpha}]=2,~[a_m,\hat{\alpha}]=0,~
[P,\hat{\alpha}]=0,~
[Q,\hat{\alpha}]=\pi i.
\end{eqnarray}
Then the following statement holds.
\begin{prop}~~~
For $c=1$, we have the free field realization
of the currents $E(z)$ and $F(z)$.
\begin{eqnarray}
E(z)&=&\epsilon(q)
:\exp\left(-\sum_{m\neq 0}
\frac{1}{[m]_q}\alpha_m z^{-m}
\right):e^{\hat{\alpha}} z^{h/2+1/2}
e^{-Q}z^{-P/r^*}
,\\
F(z)&=&\epsilon(q):
\exp\left(\sum_{m\neq 0}\frac{1}{[m]_q}\beta_m z^{-m}
\right):e^{-\hat{\alpha}} z^{-h/2+1/2}
z^{P/r+h/2r}
.
\end{eqnarray}
Here we have set
\begin{eqnarray}
\epsilon(q)=(q^{1/2}+q^{-1/2})^{-1/2}.
\end{eqnarray}
\end{prop}
Together with free field realizations of
$K(z)$ (\ref{real:EC3}), we get a free field realization of
the level one elliptic algebra $U_{q,p}(A_2^{(2)})$.

Now substituting 
the free field realization of
$E(z),\ F(z),\ K(z)$ into \eqref{real:EC1}-
\eqref{real:EC3}, we
obtain a realization of 
the half currents 
and
the $L$-operator $\hL^+(u)$ satisfying the $RLL$-relation 
\eqref{RLL} for $c=1$.
Using this $L$-operator in the "intertwining relations",
\eqref{rel:Type-I2},
\eqref{rel:Type-I5},
\eqref{rel:Type-I7},
\eqref{rel:Type-I8},
\eqref{rel:Type-I11},
 for type I and 
\eqref{rel:Type-II1}-\eqref{rel:Type-II5} for the 
 type II, one can solve them for the vertex operators.
The results are stated as follows.

\begin{thm}~~~
The highest components of the type-I and type-II
vertex operators 
$\Phi_{-}(u)$ and
$\Psi_{-}^*(u)$
are realized in terms of the free field by
\begin{eqnarray}
\Phi_{-}(z)&=&:\exp\left(-\sum_{m \neq 0}
\frac{1}{[2m]_q-[m]_q}\beta_m z^{-m}\right):
e^{\hat{\alpha}} z^{h/2+1/2}z^{-P/r-h/2r-1/r},
\label{real:Type-I1}\\
\Psi_{-}^*(z)&=&:\exp\left(
\sum_{m \neq 0}
\frac{[rm]_q}{[2m]_q-[m]_q}
\al_m z^{-m}\right):
e^{-\hat{\alpha}}z^{-h/2+1/2}e^Q z^{P/r^*+1/r^*}.
\label{real:Type-II1}
\end{eqnarray}
\end{thm}

For the other components of type-I vertex operator
$\Phi_j(u) (j=\pm,0)$,
we get the following,
by using
(\ref{rel:Type-I7}) and (\ref{rel:Type-I11}).
\begin{eqnarray}
\Phi_0(u)&=&a_{0,-}\oint_{C_0} \frac{dz'}{2\pi iz'}
\Phi_{-}(u)F(z')
\frac{[u-u'+P+\frac{h}{2}]_+}{[u-u'+\frac{1}{2}]
[
P+\frac{h}{2}+\frac{1}{2}]_+}\nonumber
\\
&=&-a_{0,-}\oint_{C_0} \frac{dz'}{2\pi iz'}
F(z')\Phi_{-}(u)
\frac{[u-u'+P+\frac{h}{2}]_+}{[u-u'-\frac{1}{2}]
[
P+\frac{h}{2}+\frac{1}{2}]_+}.
\label{real:Type-I2}
\end{eqnarray}
Here the contour $C_0$ is specified by the condition.
\begin{eqnarray}
C_0:~~|q^{-1}z|<|z'|<|p^{-1}q^{-1}z|.\nn
\end{eqnarray}
The component $\Phi_+(u)$ is given by
\begin{eqnarray}
\Phi_+(u)&=&-a_{+,-}
\oint \oint_{C_+} \frac{dz'}{2\pi iz'}
\frac{dz''}{2\pi iz''}
\Phi_-(u)F(z')F(z'')
\frac{[P+\frac{h}{2}]_+}{[
P+\frac{h}{2}-\frac{1}{2}]_+[P+\frac{h}{2}-1]_+[2P+h]}
\nonumber
\\
&\times&
\frac{[u-u'+2P+h-\frac{3}{2}][u'-u''+P+\frac{h}{2}]_+}
{[u-u'+\frac{1}{2}][u'-u''+\frac{1}{2}]}
\nonumber
\\
&=&-a_{+,-}
\oint \oint_{C^+} \frac{dz'}{2\pi iz'}
\frac{dz''}{2\pi iz''}
F(z')F(z'')
\Phi_-(u)
\frac{[P+\frac{h}{2}]_+}{[
P+\frac{h}{2}-\frac{1}{2}]_+[P+\frac{h}{2}-1]_+[2P+h]}
\nonumber
\\
&\times&
\frac{[u-u'+2P+h-\frac{3}{2}][u'-u''+P+\frac{h}{2}]_+
[u-u''+\frac{1}{2}]}
{[u-u'-\frac{1}{2}][u-u''-\frac{1}{2}][u'-u''+\frac{1}{2}]}.
\label{real:Type-I3}
\end{eqnarray}
The contour $C_+$ is specified by 
\begin{eqnarray}
C_+:~~|q^{-1}z|<|z'|<|p^{-1}q^{-1}z|,~~
|q^{-1}z|<|z''|<|p^{-1}q^{-1}z|,~~
|pqz'|<|z''|<|qz'|.\nn
\end{eqnarray}
Similarly, for Type-II vertex operators,
the component $\Psi_{0}^*(u)$ is given by
\begin{eqnarray}
\Psi_{0}^*(u)&=&a_{-,0}^*
\oint_{C_0^*} \frac{dz'}{2\pi iz'}
\Psi_-^*(u)E(z')\frac{[u-u'-P]_+^*}{[u-u'-\frac{1}{2}]^*
[P-\frac{1}{2}]_+^*}\nonumber
\\
&=&
-a_{-,0}^*
\oint_{C_0^*} \frac{dz'}{2\pi iz'}
E(z')\Psi_-^*(u)\frac{[u-u'-P]_+^*}{[u-u'+\frac{1}{2}]^*
[P-\frac{1}{2}]_+^*}.\label{real:Type-II2}
\end{eqnarray}
The contour $C_0^*$ satisfies 
\begin{eqnarray}
C_0^*:~~|q^{-1}z|<|z'|<|qz|.\nn
\end{eqnarray}
The component $\Psi_{0}^*(u)$ is given by
\begin{eqnarray}
\Psi_{+}^*(u)&=&-a_{-,+}^*
\oint \oint_{C^*_+} \frac{dz'}{2\pi iz'}
\frac{dz''}{2\pi iz''}
\Psi_-^*(u)E(z')E(z'')\frac{1}{[P-\frac{1}{2}]_+^*
[2P-2]}\nonumber\\
&\times&
\frac{[u-u'-2P+\frac{3}{2}]^*[u'-u''-P]_+^*}
{[u-u'-\frac{1}{2}]^*[u'-u''-\frac{1}{2}]^*}\nonumber
\\
&=&-
a_{-,+}^*
\oint \oint_{C^*_+} \frac{dz'}{2\pi iz'}
\frac{dz''}{2\pi iz''}
E(z')E(z'')\Psi_-^*(u)\frac{1}{[P-\frac{1}{2}]_+^*
[2P-2]}\nonumber\\
&\times&
\frac{[u-u'-2P+\frac{3}{2}]^*[u'-u''-P]_+^*[
u-u''-\frac{1}{2}]^*}
{[u-u'+\frac{1}{2}]^*
[u-u''+\frac{1}{2}]^*
[u'-u''-\frac{1}{2}]^*}.\label{real:Type-II3}
\end{eqnarray}
Here the contour $C_+^*$ is specified by the condition
\begin{eqnarray}
C_+^*:~~|q^{-1}z|<|z'|<|qz|,~~
|q^{-1}z|<|z''|<|qz|,~~|q^{-1}z'|<|z''|<|qz'|.\nn
\end{eqnarray}

~\\
{\it Remark} ~~~
The free field realizations of the vertex operators
(\ref{real:Type-I1})
-(\ref{real:Type-II3})
are the same as those of
the dilute $A_L$ model obtained in \cite{HJKOS},
up to a gauge transformation.

~\\

In addition we can verify the following commutation relation.
\begin{prop}~~~
The highest components $\Phi_-(u)$ and
$\Psi_-^*(u)$ satisfy
\begin{eqnarray}
\Phi_-(u_1)\Psi_-^*(u_2)=
\chi(u_1-u_2)
\Psi_-^*(u_2)\Phi_-(u_1).
\end{eqnarray}
Here we have set
\begin{eqnarray}
\chi(u)=-z^{-1}\frac{
\Theta_{q^6}(qz)
\Theta_{q^6}(q^2z)
}{
\Theta_{q^6}(q/z)
\Theta_{q^6}(q^2/z)
}.\label{def:chi}
\end{eqnarray}
\end{prop}


\subsection{Commutation Relations of the Vertex Operators}


We next study the commutation relations 
of the vertex operators and show that 
our realization satisfies the full 
intertwining relations for $c=1$.

\begin{thm}\label{thm:Vertexcom}~~~~
The free field realizations of
the type-I vertex operator
$\Phi_\mu(u)$ 
(\ref{real:Type-I1}),
(\ref{real:Type-I2}), 
(\ref{real:Type-I3}), and
the type-II vertex operator $\Psi_{\mu}^*(u)$
(\ref{real:Type-II1}),
(\ref{real:Type-II2}), 
(\ref{real:Type-II3}),
satisfy the following commutation relations. 
\begin{eqnarray}
&&{\Phi}^{}_{j_2}(u_2)
{\Phi}^{}_{j_1}(u_1)=
\sum_{j_1',j_2'=\pm,0}
{R}_{j_1j_2}^{j_1'j_2'}(u_1-u_2,P+h)\ 
{\Phi}^{}_{j_1'}(u_1)
{\Phi}^{}_{j_2'}(u_2)
\label{Com:Type-I},\\
&&{\Psi}^{*}_{j_1}(u_1)
{\Psi}^{*}_{j_2}(u_2)
=\sum_{j_1',j_2'=\pm,0}
{\Psi}^{*}_{j_2'}(u_2)
{\Psi}^{*}_{j_1'}(u_1)\ 
{R}^{*j_1j_2}_{j_1'j_2'}(u_1-u_2,P)
\label{Com:Type-II},\\
&&{\Phi}^{}_{j}(u_1){\Psi}^{*}_{k}(u_2)
={\chi}(u_1-u_2)\ {\Psi}^{*}_{k}(u_2)
{\Phi}^{}_{j}(u_1).\label{Com:Type-I,II}
\end{eqnarray}
Here we set
\begin{eqnarray}
{R}(u,P+h)={\mu(u)}
\bar{R}(v,P+h),~~\quad
{R}^*(u,P)
={\mu^*(u)}
\bar{R}^{*}(u,P),
\end{eqnarray}
with
\begin{eqnarray}
&&\mu(u)=z^{\frac{1}{r}-1}
\frac{\{pq^4z\}
\{pq^3z\}
\{q^3z\}
\{q^2z\}
\{pq/z\}
\{p/z\}
\{q^6/z\}
\{q^5/z\}
}{
\{pq^4/z\}
\{pq^3/z\}
\{q^3/z\}
\{q^2/z\}
\{pqz\}
\{pz\}
\{q^6z\}
\{q^5z\}
}.
\end{eqnarray}
and $\mu^*(u)=\mu(u)|_{r\to r^*}$. 
Here $\chi(u)$ is given by
(\ref{def:chi}).
\end{thm}
The proof is the similar as those of Theorem \ref{thm:HC}.

Now let us investigate the intertwining relation 
for level $c=1$.
For this purpose, we construct
a $L$-operator as a composition 
of type I and II vertex operators \cite{Miki}.

\begin{thm}~~~
For $c=1$, the components of the $L$-operator $\widehat{L}^+(u)$ \eqref{def:L}
is given by the following  product of the type-I and type-II
vertex operators.
\begin{eqnarray}
{L}_{j,k}^+(u)=g^{-1} \Psi_k^*(u+r)
\Phi_j(u+r+1/2),~~(j,k=\pm,0).\label{Mikidecomp}
\end{eqnarray}
Here we set
\begin{eqnarray}
g=-\frac{(pq^6;q^6)_\infty
(pq^5;q^6)_\infty
}{(pq^3;q^6)_\infty
(pq^2;q^6)_\infty}
\left(\frac{\{q^2p\}
\{q^3p\}
\{q^3p\}
\{q^4p\}
}{
\{p\}
\{qp\}
\{q^5p\}
\{q^6p\}}
\times
(p 
\leftrightarrow
p^*)^{-1}\right).
\end{eqnarray}
\end{thm}
The proof is similar to the one of Theorem 6.5 in \cite{KoKo}.

\noindent
{\it Remark}~~~ 
By using 
the commutation relations of the vertex operators
(\ref{Com:Type-I})-(\ref{Com:Type-I,II}) and the formula
\begin{eqnarray}
\frac{\rho^+(u)}{
\rho^{+*}(u)
}=\frac{\mu(u)\chi(\frac{1}{2}-u)}{\mu^*(u)
\chi(\frac{1}{2}+u)},
\end{eqnarray}
one can prove the $RLL=LLR^*$ relation \eqref{RLL} for $c=1$ directly.

In the same way, one can verify the ``intertwining
relations" (\ref{def:Type-I}) and (\ref{def:Type-II}) 
of vertex operators. 
\begin{cor}
~~~For 
$c=1$, the type-I and the type II vertex operators 
$\widehat{\Phi}_V(u)$,
$\widehat{\Psi}^*_V(u)$
satisfy the full intertwining relations 
\eqref{def:Type-I} and 
\eqref{def:Type-II} 
with $V=W\cong \C^3$.
\end{cor}

\section{Discussion}

Extending the construction of the elliptic algebra to the twisted 
affine Lie algebra case, we have derived the elliptic 
algebra $U_{q,p}(\att),\ p=q^{2r}$ and shown that it provides the Drinfeld 
realization of the elliptic quantum group $\Bqla(\att)$. 
Based on this, we have derived the type I and II 
vertex operators of $U_{q,p}(\att)$ and identified them with the 
vertex operators in the dilute $A_L$ model with $r=2\frac{L+1}{L+2}$.
Our result thus gives a representation theoretical foundation to the work 
\cite{HJKOS}.     

There are some open problems.
\begin{itemize}
\item[(i)] We here studied 
the simplest twisted elliptic quantum group $\Bqla(\att)$ and 
associated elliptic algebra $U_{q,p}(A_{2}^{(2)})$. 
To generalize our consideration to the higher rank cases 
associated with $A_{2n}^{(2)}$ and $A_{2n+1}^{(2)}$, or moreover to 
other types of affine Lie algebras,  
is an interesting problem. 

\item[(ii)] Our realization of the elliptic algebra 
$U_{q,p}(\att)$ based on the Drinfeld currents of 
$U_q(A_2^{(2)})$ and the Heisenberg algebra $\C\{{\cal H}\}$ 
is valid for a generic level $c$. In order to perform an algebraic analysis of 
the solvable lattice models, a free field 
realization is useful. For example, to consider a fusion 
 of the  dilute $A_L$ model, i.e. a higher spin extension, 
 we need a free field realization (Wakimoto construction)
of the elliptic algebra $U_{q,p}(A_2^{(2)})$ in higher level.

\item[(iii)] The Wakimoto realization of the affine quantum group 
$U_q(\att)$ itself is interesting. It should  be used to solve the 
$q$-KZ equation as well as the $q$-difference equation 
for the twistor $F(r,s)$, which we have solved partly 
(see Appendix \ref{Twistor}). The same thing is true for the 
other types of affine Lie algebra and should lead us to 
a proof of the conjecture on the connection
matrix of the $q$-KZ equation given by Frenkel and Reshetikhin\cite{FR}.

\item[(iv)]
It is known in some cases that the generating functions of 
the $q$-deformed Virasoro or $W$- algebras can be obtained from 
a fusion of the vertex operators of corresponding elliptic algebra   
$U_{q,p}(\g)$ \cite{JKM,JS,HJKOS,KK}. 
It is interesting to examine the same procedure in various $U_{q,p}(\g)$ 
and derive corresponding $q$-$W$ algebras. The results should be 
compared with those in \cite{EFR}.

\item[(v)] It is also an interesting problem to 
investigate the scaling limit of the half currents and 
the $L$-operators of $U_{q,p}(A_2^{(2)})$ and derive 
the vertex operators \cite{JKM,KLP}. The result should be applied to the 
Izergin-Korepin model\cite{IK} in the massless regime where  
a generic form of the correlation functions was studied in 
\cite{Koj}. The type-II vertex operators 
should provide the Zamolodchikov-Faddeev algebra
for the $A_2^{(2)}$ Toda field theory with imaginary coupling constant,
and enable us to derive the soliton $S$-matrix.
\end{itemize}
We hope to report on some of the issues listed here in the near future.


~\\
{\bf Acknowledgements}~~
The authors would like to thank Robert Weston,
Marco Rossi and colleagues in 
the Department of Mathematics,
Heriot-Watt University, where the most part of this work was done for 
their kind hospitality. They are also grateful to Michio Jimbo for sending 
his note on a construction of the half currents in $U_q(\att)$ and useful 
discussions. H.K and T.K are suported by the Grant-in-Aid for Scientific 
Research ({\bf C}) (15540033) and  the Grant-in-Aid for Young Scientist
({\bf B}) (14740107), respectively, from Japan Society for the Promotion of 
Science. H.K was also suported by JSPS-Royal Society Exchange Fellowship.


\begin{appendix}

\section{Finite Dimensional Representation}

\label{Evaluation}

The evaluation module $(\pi_w,V_w)$
in terms of the Drinfeld generators,
is defined by the following formulae.
\begin{eqnarray}
\pi_w(h)&=&2(E_{++}-E_{--}),~~\pi_w(c)=0,\\
\pi_w(a_m)&=&\frac{[m]_q}{m}(w/q)^m\left(
q^{-m}E_{++}+(1-q^m)E_{00}-q^{2m} E_{--}\right),\\
\pi_w(x_k^+)&=&(w/q)^k(aE_{+0}+q^k b E_{0-}),\\
\pi_w(x_k^-)&=&(w/q)^k(q^k b^{-1} E_{-0}+
a^{-1} E_{0+}).
\end{eqnarray}
Here we have used
$E_{++}=\left(\begin{array}{ccc}
1&0&0\\
0&0&0\\
0&0&0
\end{array}\right)$,
$E_{00}=\left(\begin{array}{ccc}
0&0&0\\
0&1&0\\
0&0&0
\end{array}\right)$ and
$E_{--}=\left(\begin{array}{ccc}
0&0&0\\
0&0&0\\
0&0&1
\end{array}\right)$.
In what follows we set
$a=b=1$.
We have
\begin{eqnarray}
&&\pi_w(x^+(z))=
E_{+0} \delta(w/qz)
+
E_{0-}
\delta(w/z),\\
&&\pi_w(x^-(z))=E_{-0}\delta(w/z)+
E_{0+} \delta(w/qz).
\end{eqnarray}
\begin{eqnarray}
\pi_w(u^+(z,p))&=&\frac{(pq^3z/w;p)_\infty}
{(pqz/w;p)_\infty}E_{++}\nonumber
\\
&+&
\frac{(pq^2z/w;p)_\infty (pq^{-1}z/w;p)_\infty}
{(pqz/w;p)_\infty (pz/w;p)_\infty}E_{00}+
\frac{(pq^{-2}z/w;p)_\infty}{(pz/w;p)_\infty}E_{--},\\
\pi_w(u^-(z,p))&=&
\frac{(pq^{-3}w/z;p)_\infty}
{(pq^{-1}w/z;p)_\infty}E_{++}\nonumber\\
&+&
\frac{(pqw/z;p)_\infty (pq^{-2}w/z;p)_\infty}
{(pw/z;p)_\infty (pq^{-1}w/z;p)_\infty}E_{00}+
\frac{(pq^2w/z;p)_\infty}{(pw/z;p)_\infty}E_{--}.
\end{eqnarray}
Let us calculate finite dimensional representation of
the elliptic current.
\begin{eqnarray}
\pi_w(e(z,p))&=&\frac{(pq^3z/w;p)_\infty
}{(pqz/w;p)_\infty}E_{+0} \delta(w/qz)+
\frac{(pq^2z/w;p)_\infty (pq^{-1}z/w;p)_\infty}
{(pqz/w;p)_\infty 
(pz/w;p)_\infty} E_{0-} \delta(w/z),\nn\\
&&\\
\pi_w(f(z,p))&=&
\frac{(pqw/z;p)_\infty
(pq^{-2}w/z;p)_\infty}{(pw/z;p)_\infty
(pq^{-1}w/z;p)_\infty}E_{-0}\delta(w/z)+
\frac{(pq^{-3}w/z;p)_\infty}{(pq^{-1}w/z;p)_\infty}
E_{0+}\delta(w/qz),\nonumber\\
\\
\pi_w(k(z,p))&=&\rho^+(q^{-r+2}z/w)\nonumber\\
&\times&\left(
E_{++}+\frac{\Theta_p(q^{r}z/w)}{\Theta_p(q^{r+2}z/w)}E_{00}+
\frac{\Theta_p(q^{r}z/w)\Theta_p(q^{r-1}z/w)}
{\Theta_p(q^{r+2}z/w)\Theta_p(q^{r+1}z/w)}E_{--}
\right),\\
\pi_w(\psi(z,p))&=&\frac{\Theta_p(q^{r+3}z/w)}
{\Theta_p(q^{r+1}z/w)}E_{++}\nonumber\\
&+&\frac{\Theta_p(q^{r+2}z/w)\Theta_p(q^{r-1}z/w)}{
\Theta_p(q^{r+1}z/w)\Theta_p(q^{r}z/w)}E_{00}+
\frac{\Theta_p(q^{r-2}z/w)}{\Theta_p(q^{r}z/w)}E_{--}.
\end{eqnarray}

\section{Twistor}

\label{Twistor}

We here consider
the difference equations
of the twistor $F(\lambda)$ for $\Bqla(\att)$. The general framework was 
given  in \cite{JKOS1}.
Let us consider the case of the affine algebra $A_2^{(2)}$.
Taking a basis $\{c,d,\al_1^{\rv}\}$ of the Cartan subalgebra
$\h$ of $A_2^{(2)}$.
We parametrize the dynamical variable $\lambda$ as
\begin{eqnarray}
\lambda-\rho=r^*d+s'c+\frac{1}{2}\left(s+\frac{r\tau}{2}\right)
\al_1^{\rv},\quad~~~~\left(\ r^*=r-c,\ \tau=-\frac{2\pi i}{\log q^{2r}}\ \right),
\end{eqnarray}
where  $\rho=3d+\frac{1}{4}\al_1^{\rv}$ is the Weyl vector.
Let us set
\begin{eqnarray}
&&\cR(z)=\Ad(z^d \otimes 1)(\cR),\\
&&F(z,p,w)=\Ad(z^d \otimes 1)(F(\lambda)),\\
&&\cR(z;p,w)=\Ad(z^d \otimes 1)(\cR(\lambda))=
\sigma(F(z^{-1};p,w))\cR(z)F(z;p,w)^{-1},
\end{eqnarray}
where $w=q^{2(s+\frac{r\tau}{2})}$.
In particular, for $z=0$, $q^{c\otimes d +d \times c}\cR(0)$
reduces to the universal $R$ matrix of 
$U_q(A_1)$.
From \cite{JKOS1}, we have the difference equation for the twistor. 
\begin{eqnarray}
&&F(pq^{2c^{(1)}}z;p,w)=
(\bar{\varphi}_w^{-1} \otimes id)(F(z;p,w))
q^T \cR(pq^{2c^{(1)}}z),
\\
&&F(0;p,w)=F_{A_1}(w),
\end{eqnarray}
where $\bar{\varphi}_w=\Ad(q^{\al_1^{\rv 2}/4}w^{\al_1^{\rv}/2})$ and
$T=\frac{1}{2}c\otimes d +\frac{1}{2}d \otimes c+\frac{1}{4}\al_1^{\rv}\otimes \al_1^{\rv}$.

We are interested in  
the vector representation $(\pi_z,V), V=\C^3$
given in Appendix \ref{Evaluation}.
We set
\begin{eqnarray}
&&F_{VV}(z;p,w)=(\pi_1 \otimes \pi_1)F(z;p,w)=
(\pi_{z_1}\otimes \pi_{z_2})(F(\lambda))
,\\
&&R_{VV}(z;p,w)=(\pi_1 \otimes \pi_1)\cR(z;p,w)=
(\pi_{z_1}\otimes \pi_{z_2})(\cR(\lambda)),\\
&&R_{VV}(z)=(\pi_1 \otimes \pi_1)\cR(z)
=(\pi_{z_1} \otimes \pi_{z_2})\cR
\end{eqnarray}
where $z=z_1/z_2$.
The trigonometric $R$-matrix $R_{VV}(z)$ 
is given as follows.
\begin{eqnarray}
&&R_{VV}(z)=\rho_{VV}(z)\bar{R}_{VV}(z),\\
&&\bar{R}_{VV}(z)=\left(\begin{array}{ccccccccc}
1&0&0&0&0&0&0&0&0\\
0&b(z)&0&c(z)&0&0&0&0&0\\
0&0&d(z)&0&e(z)&0&f(z)&0&0\\
0&z\ c(z)&0&b(z)&0&0&0&0&0\\
0&0&-q^2z\ e(z)&0&j(z)&0&e(z)&0&0\\
0&0&0&0&0&b(z)&0&c(z)&0\\
0&0&z\ n(z)&0&-q^2z\ e(z)&0&d(z)&0&0\\
0&0&0&0&0&z\ c(z)&0&b(z)&0\\
0&0&0&0&0&0&0&0&1
\end{array}\right),\nonumber\\
&&b(z)=-\frac{q(1-z)}{1-q^2z},\qquad
c(z)=\frac{1-q^2}{1-q^2z},\qquad
d(z)=\frac{(1-z)q^2(1-qz)}{(1-q^2z)(1-q^3z)},\nn\\
&&e(z)=\frac{i(1-q^2)q^{\frac{1}{2}}(1-z)}
{(1-q^2z)(1-q^3z)},\qquad
f(z)=
\frac{(1-q^2)(1+q-q^3z-qz)}{
(1-q^2z)(1-q^3z)},\nn\\
&&j(z)=-\frac{q(1-z)}{1-q^2z}+
\frac{(1-q^2)(1-q^3)z}{(1-q^2z)(1-q^3z)},\qquad
n(z)=\frac{(1-q^2)(1+q^2-q^3z-q^2z)}{
(1-q^2z)(1-q^3z)}.\nn
\end{eqnarray}
The function $\rho_{VV}(z)$
is given by
\begin{eqnarray}
\rho_{VV}(z)&=&q^{-1}\frac{(1/z;q^6)_\infty
(q/z;q^6)_\infty
(q^5/z;q^6)_\infty(q^6/z;q^6)_\infty
}{
(q^2/z;q^6)_\infty
(q^3/z;q^6)^2_\infty
(q^4/z;q^6)_\infty
}.
\end{eqnarray}
Noting $\pi_1(c)=0$, 
we have the difference equation 
\begin{eqnarray}
F_{VV}(pz;p,w)^t=R_{VV}(pz)^t
K(D_w \otimes 1)^{-1}
F_{VV}(z;p,w)^t (D_w \otimes 1),\label{diff:Twistor}
\end{eqnarray}
where $X^t$ means the transpose of $X$,
and we have set
\begin{eqnarray}
K&=&\Diag(q,1,q^{-1},1,1,1,q^{-1},1,q),\\
D_w&=&\Diag(q^{-1}w^{-1},1,q^{-1}w).
\end{eqnarray}
From the from of $R_{VV}(z)$, one can set 
\bea
&&F(z;p,w)=f(z)\left(\begin{array}{ccccccccc}
1&0&0&0&0&0&0&0&0\\
0&X^{(+)}_{11}(z)&0&X^{(+)}_{12}(z)&0&0&0&0&0\\
0&0&Y_{11}(z)&0&Y_{12}(z)&0&Y_{13}(z)&0&0\\
0&X^{(+)}_{21}(z)&0&X^{(+)}_{22}(z)&0&0&0&0&0\\
0&0&Y_{21}(z)&0&Y_{22}(z)&0&Y_{23}(z)&0&0\\
0&0&0&0&0&X^{(-)}_{11}(z)&0&X^{(-)}_{12}(z)&0\\
0&0&Y_{31}(z)&0&Y_{32}(z)&0&Y_{33}(z)&0&0\\
0&0&0&0&0&X^{(-)}_{21}(z)&0&X^{(-)}_{22}(z)&0\\
0&0&0&0&0&0&0&0&1
\end{array}\right).\nn
\ena
Then the $q$-difference equation (\ref{diff:Twistor})
is equivalent to the following 
equations.
\bea
&&f(pz)=q\rho_{VV}(pz)f(z),\\
&&\nn\\
&&\left(\matrix{X^{(\pm)}_{11}(pz)&X^{(\pm)}_{12}(pz)\cr
                X^{(\pm)}_{21}(pz)&X^{(\pm)}_{22}(pz)\cr}\right)=
 q^{-1}\left(\matrix{X^{(\pm)}_{11}(z)& q^{\pm 1}w^{}X^{(\pm)}_{12}(z)\cr
                q^{\mp1}w^{-1}X^{(\pm)}_{21}(z)&X^{(\pm)}_{22}(pz)\cr}\right)
\left(\matrix{b(pz)&c(pz)\cr
               pz\ c(pz)&b(pz)\cr}\right),\nn\\
&&\\
&&\left(\matrix{Y_{11}(pz)&Y_{12}(pz)&Y_{13}(pz)\cr
                Y_{21}(pz)&Y_{22}(pz)&Y_{23}(pz)\cr
                Y_{31}(pz)&Y_{32}(pz)&Y_{33}(pz)\cr}\right)\nn\\
&&\qquad =q^{-2}
\left(\matrix{Y_{11}(z)&wY_{12}(z)&w^2Y_{13}(pz)\cr
                qw^{-1}Y_{21}(z)&qY_{22}(z)&qwY_{23}(z)\cr
                w^{-2}Y_{31}(z)&w^{-1}Y_{32}(z)&Y_{33}(z)\cr}\right)           
\left(\matrix{d(pz)&e(pz)&f(pz)\cr
                -q^2pz\ e(pz)&j(pz)&e(pz)\cr
                pz\ n(pz)&-q^2pz\ e(pz)&d(pz)\cr}\right).\nn\\
\lb{3t3}\ena
The two $2\times 2$ matrix equations for $X^{(\pm)}_{ij}(z)$ 
are the same as the one appeared in  the $\slth$ case \cite{JKOS1},
if we change $b(z)$ to $-b(z)$ and make the following identification.
\be
&&q^{\pm1}w=w^{-1}_{\slt}
\quad 
{\rm i.e.}\quad -s_{\slt}=s+\frac{r\tau}{2}\pm\frac{1}{2}, 
\en
where $w_{\slt}$ and $s_{\slt}$ 
denote $w$ and $s$ in \cite{JKOS1}, respectively. Hence from the elliptic $R$ 
matrix for $\Bqla(\slth)$ ((4.18) in \cite{JKOS1}), we determine the 
following parts of our elliptic $R$-matrix
\bea
&&\left(\matrix{R^{+0}_{+0}&R^{0+}_{+0}\cr
R^{+0}_{0+}&R^{0+}_{0+}\cr}\right)
=\left(\matrix{-\frac{[s+3/2]_+[s-1/2]_+}{[s+1/2]_+^2}\frac{[u]}{[u+1]}&e^{\frac{\pi i u}{r}}\frac{[1][s+1/2-u]_+}{[s+1/2]_+[u+1]}\cr
e^{-\frac{\pi i u}{r}}\frac{[1][s+1/2+u]_+}{[s+1/2]_+[u+1]}&-\frac{[u]}{[u+1]}\cr} \right),\\
&&\nn\\
&&\left(\matrix{R^{0-}_{0-}&R^{-0}_{0-}\cr
R^{0-}_{-0}&R^{-0}_{-0}\cr}\right)
=\left(\matrix{-\frac{[s-3/2]_+[s+1/2]_+}{[s-1/2]_+^2}\frac{[u]}{[u+1]}&
e^{\frac{\pi i u}{r}}\frac{[1][s-1/2-u]_+}{[s-1/2]_+[u+1]}\cr
e^{-\frac{\pi i u}{r}}\frac{[1][s-1/2+u]_+}{[s-1/2]_+[u+1]}&-\frac{[u]}{[u+1]}\cr} \right). 
\ena
By a gauge transformation, these yields the corresponding matrix elements 
in \eqref{rmat}.

As for the $3\times 3$ part, we have no known solutions.
The Wakimoto realization of $U_{q}(\att)$ should be 
useful to solve the $q$-KZ equation for the intertwining operators 
(vertex operators) of $U_{q}(\att)$ as well as \eqref{3t3}.

\section{Proof of the Relation \eqref{rel:HC8}}
\label{Proof}

Let us set 
\bea
&&h(v)=-\frac{[v+1]^*[v-1/2]^*}{[v-1]^*[v+1/2]^*}.
\ena
In the integrand of the half current $E^+_{-,+}(u)$ 
\eqref{def:HC4}, 
we call $E(z')E(z'')$ the operator part, and the ratio of the product of 
the theta functions the coefficient part.
We keep coefficient parts in the right of operator parts.
According to the relation \eqref{def:EA4}, we have the equality
\be
\oint
{d z'\over2\pi iz'}{dz''\over2\pi iz''}
E(z')E(z'')A(u,u'')
=\oint
{dz'\over2\pi iz'}{dz''\over2\pi iz''}
E(z')E(z'')h(u''-u')A(u'',u'),
\en
when the integration contours for $z'$ and $z''$ are the same. Here we set 
$z'=q^{2u'},\ z''=q^{2u''}$.
We define `weak equality' in the following sense\cite{AJMP}. 
The two coefficient functions $A(u',u'')$ and $B(u',u'')$ coupled to
$E(z')E(z'')$ in integrals are equal in weak sense if
\be
A(u',u'')+h(u''-u')A(u'',u')
=B(u',u'')+h(u''-u')B(u'',u').
\en
We write the weak equality as
\be
A(u',u'')\sim B(u',u'').
\en
To prove the equality \eqref{rel:HC8},
it is enough to show the equalities of coefficient parts in weak sense.

Setting $z_i=q^{2u_i}\ (i=1,2)$ and $u=u_1-u_2$, let us consider {RHS}-{LHS}
 of \eqref{rel:HC8} given as follows.
\be
&&\oint \frac{dz'}{2\pi iz'}
\oint\frac{dz''}{2\pi i z''} 
E(z')E(z'')F(u_1,u_2,u',u'',L),
\en
where 
\be
&&F(u_1,u_2,u',u'',L)\nn\\
&&=\frac{[u_2-u'-2P+2+c/2]^*[u'-u''-P]_+^*[1+u]^*[u+3/2]^*[1]^{*2}}
{[u_2-u'+c/2]^*[u'-u''-1/2]^*[u]^*[2P-2]^*[P-1/2]_+^*[u+1/2]^*}\nn\\
&&+\frac{[u_2-u'-P+(c+1)/2]_+^*[u_1-u'+1+c/2]^*[u_1-u"-P+(c+1)/2]_+^*
[u+P+1]_+^*[1]^{*}}
{[u_2-u'+c/2]^*[u_1-u'+c/2]^*[u_1-u''+c/2]^*[P-1/2]_+^{*2}[P+1/2]_+^*[u+1/2]^*}\nn\\
&&-\frac{[u_1-u'-2P+2+c/2]^*[u'-u''-P]_+^*[1]^{*2}}
{[u_1-u'+c/2]^*[u'-u''-1/2]^*[P-1/2]_+^*[u+1/2]^*}\nn\\
&&\qquad\times \left(\frac{[u+2P-1]^*[1]^*[u+3/2]^*}{[u]^*[2P-1]^*[2P-2]^*}
+\frac{[P]_+^*[u+2P+1/2]^*[1]^*}{[2P]^*[2P-1]^*[P-1]_+^*}
\right)\nn\\ 
&&-\frac{[u_2-u'-2P+c/2]^*[u'-u''-P-1]_+^*[u_1-u'+1+c/2]^*[u_1-u''+1+c/2]^*[1]^{*2}}
{[u_2-u'+c/2]^*[u'-u''-1/2]^*[u_1-u'+c/2]^*[u_1-u''+c/2]^*[P+1/2]_+^*[2P]^*}.
\en
We will show that $F(u_1,u_2,u',u'',L)\sim 0$.
For this purpose, we consider the function of $u'$ defined by
\be
&&F(u')=F(u_1,u_2,u',u'',L)+h(u''-u')F(u_1,u_2,u'',u',L).
\en
Then it is not so hard to see that $F(u')$ is a quasi-periodic function having 
zeros at least at $u'=u''$ and $u'=u''+1$. The quasi-periodicity is given by
\be
&&F(u'+\tau^* r^*)=-e^{-\frac{2\pi i}{r}(P-3/2)}F(u'),\nn\\
&&F(u'+r^*)=F(u').
\en
Therefore if we set 
\be
&&G(u')=F(u')\frac{[u'-u''-P+3/2]^*}{[u'-u'']^*}, 
\en
$G(u')$ is a doubly periodic function of $u'$ and $G(u''+1)=0$. 
It is then enough to show that $G(u')$ is an entire function.

In $G(u')$, some of terms have the first order poles at 
$u'=u_1+c/2,\ u_2+c/2, u''+1/2, u''-1$. We checked that all 
the residues of the function $G(u')$ at these poles vanish. 
For example, at $u'=u_2+c/2$ 
the residue is given by
\be
{\rm Res}_{u'=u_2+c/2}G(u')\frac{dz'}{2\pi i z'}
&=&-\frac{[u+1]^*[u_2-u''-P+c/2]_+^*[u+3/2]^*[1]^{*2}}
{[u]^*[u_2-u''-1/2+c/2]^*[P-1/2]_+^*[u+1/2]^*}\nn\\
&& +\frac{[u+1]^*[u_1-u''-P+(c+1)/2]_+^*[u+P+1]_+^*[1]^{*3}}
{[u]^*[u_1-u''+c/2]^*[P+1/2]_+^*[P-1/2]_+^*[u+1/2]^*}\nn\\
&&+\frac{[u+1]^*[u_1-u''+1+c/2]^*[u_2-u''-P-1+c/2]_+^*[1]^{*2}}
{[u]^*[u_1-u''+c/2]^*[u_2-u''-1/2+c/2]^*[P+1/2]_+^*}.
\en
One can apply the following 
theta function identity to combine the 1st and the 3rd terms.
\be
&&[u+x]^*[u-x]^*[v+y]_+^*[v-y]_+^*-[u+y]^*[u-y]^*[v+x]_+^*[v-x]_+^*\nn\\
&&\qquad =-[x-y]^*[x+y]^*[u+v]_+^*[u-v]_+^*.
\en
We thus get 
\be
&&{\rm 1st}+{\rm 3rd}=-\frac{[u+1]^*[u_1-u''-P+(c+1)/2]_+^*[u+P+1]_+^*[1]^{*3}}
{[u]^*[u_1-u''+c/2]^*[P+1/2]_+^*[P-1/2]_+^*[u+1/2]^*}.
\en
Therefore ${\rm Res}_{u'=u_2+c/2}G(u')\frac{dz'}{2\pi i z'}=0$.
The other cases can be treated in the similar way.

\end{appendix}

\end{document}